\documentclass[12pt,a4paper]{article}
\usepackage[english]{babel}
\usepackage{amsmath}
\usepackage{amssymb}
\usepackage{amsfonts}
\usepackage{hyperref}

\begin{document}

%\binoppenalty=10000 \relpenalty=10000

\renewcommand{\refname}{References}
\renewcommand{\contentsname}{Contents}

\begin{center}
{\huge Solvability of a steady boundary-value problem for the equations of one-temperature viscous compressible heat-conducting bifluids}
\end{center}

\medskip

\begin{center}
{\large Alexander Mamontov,\quad Dmitriy Prokudin}
\end{center}

\medskip

\begin{center}
{\large October 18, 2017}
\end{center}

\medskip

\begin{center}
{%\large
Lavrentyev Institute of Hydrodynamics, \\ Siberian Branch of the Russian Academy of Sciences\\
pr. Lavrent'eva 15, Novosibirsk 630090, Russia}
\end{center}

\medskip

\begin{center}
{\bfseries Abstract}
\end{center}

%\medskip

\begin{center}
\begin{minipage}{110mm}
We consider a boundary-value problem describing the steady motion of a two-component mixture of viscous compressible heat-conducting fluids in a bounded domain. We make no simplifying assumptions except for postulating the coincidence of phase temperatures (which is physically justified in certain situations), that is, we retain all summands in equations that are a natural generalization of the Navier-Stokes-Fourier model of the motion of a one-component medium. We prove the existence of weak generalized solutions of the problem.
\end{minipage}
\end{center}

\bigskip


\newpage

\tableofcontents

\bigskip

\section{Introduction}

\noindent\indent The description of the motion of multicomponent media is an interesting and rather little-studied problem both in physics/mechanics and in mathematics. There
is no standard approach to simulating these motions, nor is there any developed mathematical theory concerning the existence, uniqueness and properties of solutions of initial-boundary value problems arising in this simulation. The aims of our paper do not include a detailed survey of these problems as a whole. To some extent one can form an impression of them from the monographs \cite{nigm} and \cite{raj} and also from the survey in the paper \cite{france}. Our object is to carry out a mathematical investigation, although due to the physical origin of the problem, short explanations concerning the mechanical meaning of some aspects are unavoidable, and we shall give these explanations when necessary.

In the paper, we choose one of the numerous versions of simulating the motion of two-velocity (binary) fluid mixtures, namely, a homogeneous mixture of two viscous compressible heat-conducting fluids and a two-velocity one-temperature model. This means that both constituents of the mixture are present at any point of the space and with the same phase, and each of them has its own local velocity. The interaction between the components is studied using the viscous friction and the exchange between momenta, and also using the heat exchange (in heat-conducting models). We note that, from a mathematical point of view, both this model and numerous others have received little attention (when compared to the corresponding theory for one-component media). Some details can be found in the survey \cite{france}. Some explanations are given below. In connection with the model considered in this paper (and similar models), we note that, concerning binary two-velocity mixtures, until recently the only results were for approximate models, which do not consider temperatures \cite{frehsesiam}, \cite{frehseAM}, \cite{frehseweig}. Papers have recently appeared concerning the complete model, first for the barotropic model \cite{kucher} and then for the heat-conducting (two-temperature) model \cite{smj12}. The latter paper still contains some simplification, namely, the terms (in the energy equations) corresponding to viscous friction are absent. This is connected not only with mathematical complications (which also occur in the one-component case) but also with the physical correctness of the model (for details, see \cite{france}). In the paper, we consider a one-temperature model in which the physical discrepancies (in dissipative terms) that are typical of mixtures do not occur, and thus only mathematical difficulties remain, and we have been able to overcome these. Thus, we have succeeded in obtaining a first result concerning the mathematical correctness of a complete heat-conducting model of a binary mixture in the case of multidimensional motions.  We note that there are results in the one-dimensional case concerning heat-conducting mixtures (\cite{papin}, \cite{petrov}). However, they deal with approximate models, in which the viscosity matrix is diagonal. To conclude this survey, let us mention several papers concerning barotropic flows: \cite{mamprok.arxiv1}, \cite{mamprok1d.smz1}, \cite{mamprok1d.smz2}, \cite{mamprok.arxiv2}, \cite{mamprok.arxiv3}.

When modeling the motion of mixtures, there is a risk of encumbering the notation because additional indices can occur (corresponding to indexing the components of the mixture). A significant simplification is achieved here by using invariant notation (excluding explicit mention of the components of vectors and tensors) which we use throughout the paper, refining its rules to avoid ambiguity. Thus, if $\boldsymbol{a}$ and $\boldsymbol{b}$ are vectors (``columns'') of dimension~$n$ and ${\mathbb A}$ and ${\mathbb B}$  are tensors of
the second rank (``matrices'') acting on~${\mathbb R}^n$, then
$$\boldsymbol{a}\cdot\boldsymbol{b}=\sum\limits_{i=1}^n a_i b_i,\quad {\mathbb A}:{\mathbb B}=\sum\limits_{i,j=1}^n A_{ij}B_{ij},\quad
{\rm div}\boldsymbol{a}=\sum\limits_{i=1}^n\frac{\partial a_i}{\partial x_i},$$
${\mathbb A}\boldsymbol{a}$  and ${\rm div}{\mathbb A}$  are vectors (``columns'') with components
$$({\mathbb A}\boldsymbol{a})_i=\sum\limits_{j=1}^n A_{ij}a_j,\quad ({\rm div}{\mathbb A})_i=\sum\limits_{j=1}^n \frac{\partial A_{ij}}{\partial x_j},$$
and $\boldsymbol{a}\otimes\boldsymbol{b}$  is the tensor with components $(\boldsymbol{a}\otimes\boldsymbol{b})_{ij}=a_i b_j$.

\section{Statement of the problem and the main\\ result}

\noindent\indent The model to be used here is a natural generalization of the\linebreak Navier-Stokes-Fourier model of the motion of a one-component viscous compressible heat-conducting fluid, and has revived significant mathematical study over the last two decades (some results of this work can be found in the monographs 
\cite{lions98}, \cite{feir03}, \cite{novstrs04}, \cite{feirnov09} and \cite{plotsok12}. To formulate the model of the motion of a mixture as a generalization of the one-component model, significant efforts are needed, and these are reflected, for example, in the monographs \cite{nigm} and \cite{raj}. The two-velocity model arising in this way admits certain variations. In particular, there remains an arbitrariness in modeling the effects of temperature. One can assume either a distinction between the temperatures in different components of the mixture or the coincidence of these temperatures. The latter assumption is used in the paper, and it is justifiable under certain physical conditions (see, for example, \cite{voinov}, \cite{jum} and \cite{gard}).

Omitting intermediate considerations related to the construction of our model, we immediately formulate the mathematical problem to be considered. Let a mixture of two viscous compressible heat-conducting fluids occupy a bounded domain $\Omega\subset{\mathbb R}^3$. We regard the assumptions on the number of components of the mixture, on the fact that the flow is three-dimensional, and on the stationarity of the flow as inessential (the corresponding generalizations of the result involve no fundamental difficulties). We need the following physical quantities described by five (or with regard to the dimensions of the vectors, by nine) functions defined on $\Omega$:  the scalar density fields $\rho_i\geqslant 0$ and the vector velocity fields $\boldsymbol{u}^{(i)}$ for every component of the
mixture, $i=1,2$,  and also the scalar field of the temperature of the mixture, $\theta>0$. To find these quantities, one must solve two continuity equations,
\begin{align}\label{continuity}{\rm div}(\rho_i \boldsymbol{u}^{(i)})=0\quad\text{on}\quad\Omega,\quad i=1,2,\end{align}
two vector (that is, six scalar) equations for momenta,
\begin{align}\label{momentum}\sum\limits_{j=1}^2 L_{ij}\boldsymbol{u}^{(j)}+{\rm div}(\rho_i \boldsymbol{u}^{(i)}\otimes \boldsymbol{u}^{(i)})+\nabla p_i=
\rho_i \boldsymbol{f}^{(i)}+\boldsymbol{J}^{(i)}\quad\text{on}\quad\Omega,\; i=1,2,\end{align}
and one equation for the total energy of the mixture,
\begin{align}\label{energy}\begin{array}{c}\displaystyle
{\rm div}\left(\sum\limits_{i=1}^2\rho_i E_i \boldsymbol{u}^{(i)}\right)+{\rm div}\left(\sum\limits_{i=1}^2 p_i \boldsymbol{u}^{(i)}\right)-
{\rm div}\left(\sum\limits_{i=1}^2 {\mathbb P}^{(i)} \boldsymbol{u}^{(i)}\right)=\\ \\
\displaystyle =\sum\limits_{i=1}^2 \rho_i \boldsymbol{f}^{(i)}\cdot\boldsymbol{u}^{(i)}-2{\rm div}\boldsymbol{q}\quad\text{on}\quad\Omega.
\end{array}\end{align}
In the equations \eqref{continuity}-\eqref{energy} we use the following notation. $\boldsymbol{f}^{(i)}$ are the known external mass forces.
$p_i$ is the pressure of the $i$th component, for which the validity of the following constitutive equation is assumed:
\begin{align}\label{pressure}p_i=\rho_i^\gamma+\rho_i\theta,\quad i=1,2,\end{align}
in which the adiabatic exponent $\gamma$ is assumed to be the same for both components (which is inessential) and sufficiently large (precise conditions are given
below and can be weakened); the nature of the equation \eqref{pressure}  is fairly standard in Navier-Stokes-Fourier theory (see, for example, 
\cite{novpok11}, \cite{muchapok09} and \cite{muchapok10}). $\boldsymbol{J}^{(i)}$ is the momentum exchange (between the components of the mixture) defined in the standard way:
\begin{align}\label{momsupply}\boldsymbol{J}^{(i)}=(-1)^i a(\boldsymbol{u}^{(1)}-\boldsymbol{u}^{(2)}),\quad i=1,2,\quad a={\rm const}>0.\end{align}
$\boldsymbol{q}$ is the heat-flux vector defined by the Fourier law
\begin{align}\label{heatflux}\boldsymbol{q}=-k(\theta)\nabla\theta\end{align}
with the heat-conducting coefficient assumed to be of the form
\begin{align}\label{heatcondcoeff}k(\theta)=1+\theta^m,\end{align}
where the constant $m$  is specified below. $E_i$ is the total specific energy of the $i$th
component of the mixture, defined as the sum of the kinetic and internal energies,
\begin{align}\label{totalenergy}E_i=\frac{1}{2}|\boldsymbol{u}^{(i)}|^2+U_i,\quad i=1,2,\end{align}
where the specific internal energy $U_i$ is given by the constitutive equation
\begin{align}\label{intenergy}U_i=\frac{1}{\gamma-1}\rho_i^{\gamma-1}+\theta,\quad i=1,2.\end{align}
${\mathbb P}^{(i)}$ is the viscous part of the stress tensor of the $i$th component of the mixture
given by the constitutive equation
\begin{align}\label{stress}{\mathbb P}^{(i)}=\sum\limits_{j=1}^2\widehat{{\mathbb P}}^{(ij)},\quad\text{where}\quad
\widehat{{\mathbb P}}^{(ij)}=\lambda_{ij}{\rm div}\boldsymbol{u}^{(j)}{\mathbb I}+2\mu_{ij}{\mathbb D}(\boldsymbol{u}^{(j)}),
\quad i,j=1,2,\end{align}
in which the (constant) viscosity coefficients $\lambda_{ij}$ and $\mu_{ij}$ necessarily satisfy certain
conditions given below, ${\mathbb I}$ stands for the identity tensor, and\linebreak ${\mathbb D}(\boldsymbol{v})=((\nabla\otimes\boldsymbol{v})+(\nabla\otimes\boldsymbol{v})^T)/2$ is the strain tensor of the vector field $\boldsymbol{v}$ (the superscript $T$ stands for transposition). Finally, we use the following notation for Lam\'e operators:
\begin{align}\label{lame}L_{ij}=-(\lambda_{ij}+\mu_{ij})\nabla{\rm div}-\mu_{ij}\Delta,\quad i,j=1,2,\end{align}
and thus ${\rm div}{\mathbb P}^{(i)}=-\sum\limits_{j=1}^2 L_{ij}\boldsymbol{u}^{(j)}$, $i=1,2$. The thermodynamical constitutive equations \eqref{pressure},
\eqref{heatcondcoeff} и \eqref{intenergy} are chosen in this way for definiteness, and can be generalized (naturally, in the framework of fundamental constraints such as the Gibbs relation). The relations \eqref{momsupply} and \eqref{stress} are the prevalent version for modeling the dynamics of two-component mixtures, and express the principles of mechanical interaction of the components that are formulated above. Finally, \eqref{heatflux} and \eqref{totalenergy} represent standard physical laws.

To the equations \eqref{continuity}-\eqref{energy} one must add boundary conditions for the velocities and temperatures, for example,
\begin{align}\label{boundvelocity}\boldsymbol{u}^{(i)}=0\quad\text{on}\quad\partial\Omega,\quad i=1,2,\end{align}
that is, the boundary $\partial\Omega$ of $\Omega$  is assumed to be an impermeable rigid wall (which can however be readily generalized),
\begin{align}\label{boundtemper}2k(\theta)\nabla\theta\cdot\boldsymbol{n}+L(\theta)(\theta-\widehat{\theta})=0\quad\text{on}\quad\partial\Omega\end{align}
(that is, there is a heat exchange with the external medium which has a known temperature distribution $\widehat{\theta}>0$), and also additional conditions on the densities,
which as is customary we adopt in the form
\begin{align}\label{addconddensity}\int\limits_\Omega \rho_i \,d\boldsymbol{x}=M_i,\quad i=1,2,\end{align}
where the positive constants $M_i$ express the total masses of the components of the mixture and are assumed to be known. Under the assumptions \eqref{boundtemper}, the
symbol $\boldsymbol{n}$  stands for the outward-pointed unit normal to $\partial\Omega$, and the coefficient of the boundary heat exchange is adopted in the form
\begin{align}\label{boundheatex}L(\theta)=1+\theta^{m-1}.\end{align}

Thus, the subject of our investigation has been specified: it is the\linebreak boundary-value problem \eqref{continuity}--\eqref{boundheatex}, to be referred to in what follows as Problem ${\mathcal H}$.

It can be seen that the proposed model of the motion of a binary mixture is not the simple union of two Navier-Stokes-Fourier systems describing the motion of each of the components, because these systems interact in the higher order terms, namely, in the summands responsible for the viscous friction between the components. The viscosity coefficients form two matrices $\Lambda=\{\lambda_{ij}\}_{i,j=1}^2$ and $M=\{\mu_{ij}\}_{i,j=1}^2$ (of bulk and shear viscosity) whose off-diagonal entries are responsible for the above interaction. An important role is also played by the matrix of total viscosities $N=\Lambda+2M$ (with the entries $\nu_{ij}=\lambda_{ij}+2\mu_{ij}$). If these matrices are diagonal, then the interaction of the systems involves only the lower order terms, and the corresponding problem contains no substantial mathematical difficulties when compared with one-component motion (although it is still physically interesting). We consider the case of non-diagonal viscosity matrices, in which it is impossible to immediately transfer the approaches developed in the Navier-Stokes-Fourier theory of one-component fluids. This is related to the mathematical interest in the problem posed above, and it is this interest that stimulated the publication of the paper.

For the thermodynamical consistency of the model formulated above, the viscosity matrices must satisfy certain conditions of positive or non-negative definiteness
(see \cite{france}). Here we assume the validity of the following conditions (which are close to minimal):
\begin{align}\label{positiveviscos}M>0,\qquad 3\Lambda+2M\geqslant 0,\end{align}
ensuring the validity if the inequalities
\begin{align}\label{nonnegentr}\sum\limits_{i=1}^2{\mathbb P}^{(i)}:(\nabla\otimes\boldsymbol{u}^{(i)})\geqslant 0\end{align}
(corresponding to the non-negativeness of the entropy production), and
\begin{align}\label{coerciveviscos}\sum\limits_{i,j=1}^2 \int\limits_\Omega L_{ij}\boldsymbol{u}^{(j)}\cdot \boldsymbol{u}^{(i)}\,d\boldsymbol{x}\geqslant C_0
\sum\limits_{i=1}^2 \int\limits_\Omega|\nabla\otimes\boldsymbol{u}^{(i)}|^2 \,d\boldsymbol{x}\end{align}
(ensuring the ellipticity property, which is important from a mathematical point of view), where
$$2C_0=(\mu_{11}+\mu_{22})-\sqrt{(\mu_{11}-\mu_{22})^2+(\mu_{12}+\mu_{21})^2}.$$
We note that \eqref{positiveviscos} yields that $N>0$. Moreover, to simplify the application of the method of effective viscous flux (which is at the heart of the modern Navier-Stokes-Fourier theory) in the matrix case, we make the additional technical assumption that the matrix of total viscosities is triangular, namely,
\begin{align}\label{triangleviscos}\lambda_{12}+2\mu_{12}=0.\end{align}
This is not used in the subsequent construction except at the point indicated above (see stages 4.2 and 4.3 of the proof of Theorem 2.3 in section 5). However, the problem
of the efficiency of the method of effective viscous flux remains open, and may be of critical difficulty in the case of arbitrary viscosity matrices.

The object of the paper is the construction of a weak generalized solution of Problem ${\mathcal H}$ understood in the standard way, quite in the spirit of the theory of one-component viscous gases. To be precise, we give a rigorous definition.

{\bfseries Definition 2.1.} {\it By a weak generalized solution of Problem ${\mathcal H}$ we mean a pair of
non-negative functions $\rho_i\in L_{2\gamma}(\Omega)$, $i=1,2$,  a positive function $\theta\in W^1_2(\Omega)\bigcap L_{3m}(\Omega)\bigcap L_{2m}(\partial\Omega)$, and a pair of vector fields $\boldsymbol{u}^{(i)}\in \overset{\circ}{W^1_2}(\Omega)$,
$i=1,2$,  satisfying
the conditions \eqref{continuity}--\eqref{boundheatex}  in the following sense.
\begin{enumerate}
  \item[(${\mathcal H}1$)] The densities $\rho_i$ satisfy the continuity equations \eqref{continuity}, which means that
the integral identities
$$\int\limits_\Omega \rho_i\boldsymbol{u}^{(i)}\cdot\nabla\psi_i\,d\boldsymbol{x}=0,\qquad i=1,2,$$
hold for every $\psi_i\in C^\infty(\overline{\Omega})$, and the conditions \eqref{addconddensity} hold.
  \item[(${\mathcal H}2$)] The velocities $\boldsymbol{u}^{(i)}$ satisfy the equations for momenta \eqref{momentum} (with the constitutive equations and notation \eqref{pressure},
\eqref{momsupply}, \eqref{stress}, \eqref{lame}), which means that the integral identities
$$\sum\limits_{j=1}^2\left(\mu_{ij}\int\limits_\Omega(\nabla\otimes\boldsymbol{u}^{(j)}):(\nabla\otimes\boldsymbol{\varphi}^{(i)})\,d\boldsymbol{x}+\right.$$
$$\left.+  (\lambda_{ij}+\mu_{ij})\int\limits_\Omega({\rm div}\boldsymbol{u}^{(j)})({\rm div}\boldsymbol{\varphi}^{(i)})\,d\boldsymbol{x}\right)-$$
$$-\int\limits_\Omega (\rho_i\boldsymbol{u}^{(i)}\otimes \boldsymbol{u}^{(i)}):(\nabla\otimes\boldsymbol{\varphi}^{(i)})\,d\boldsymbol{x}-
  \int\limits_\Omega\rho_i^\gamma {\rm div}\boldsymbol{\varphi}^{(i)}\,d\boldsymbol{x}-$$
$$-\int\limits_\Omega\rho_i\theta{\rm div}\boldsymbol{\varphi}^{(i)}\,d\boldsymbol{x}=\int\limits_\Omega(\rho_i\boldsymbol{f}^{(i)}+\boldsymbol{J}^{(i)})\cdot
  \boldsymbol{\varphi}^{(i)}\,d\boldsymbol{x},\qquad i=1,2,$$
hold for arbitrary vector fields $\boldsymbol{\varphi}^{(i)}\in C_0^\infty(\Omega)$ $($the boundary conditions \eqref{boundvelocity} hold
in the sense of the function class$)$.
  \item[(${\mathcal H}3$)] The temperature $\theta$ satisfies the energy equation \eqref{energy} $($with the constitutive equations and notation \eqref{heatflux}-\eqref{intenergy}$)$  and the boundary condition \eqref{boundtemper} $($with the constitutive equation \eqref{boundheatex}$)$, which means that the integral identity
$$-\sum\limits_{i=1}^2 \int\limits_\Omega \rho_i E_i \boldsymbol{u}^{(i)}\cdot\nabla\eta \,d\boldsymbol{x}-
   \sum\limits_{i=1}^2 \int\limits_\Omega p_i \boldsymbol{u}^{(i)}\cdot\nabla\eta \,d\boldsymbol{x}+$$
$$+\sum\limits_{i=1}^2 \int\limits_\Omega {\mathbb P}^{(i)}:(\boldsymbol{u}^{(i)}\otimes\nabla\eta)\,d\boldsymbol{x}=
   \sum\limits_{i=1}^2 \int\limits_\Omega \rho_i \boldsymbol{f}^{(i)}\cdot\boldsymbol{u}^{(i)}\eta \,d\boldsymbol{x}+$$
$$+2\int\limits_\Omega\boldsymbol{q}\cdot\nabla\eta \,d\boldsymbol{x}-\int\limits_{\partial\Omega} L(\theta)(\theta-\widehat{\theta})\eta d\sigma$$
holds for every $\eta\in C^\infty(\overline{\Omega})$.
\end{enumerate}}

{\bfseries Remark 2.2.} {\it It is known from the theory of transport equations and the Navier-Stokes-Fourier theory $($see, for example, \cite{diperna}, \cite{lions98}, \cite{novstrs04}$)$ that all weak solutions of the continuity equations $($in the sense of part ${\mathcal H}1$  of Definition 2.1$)$ are automatically so-called renormalized solutions, that is, they satisfy the renormalized equations that are formally obtained from\eqref{continuity} by multiplying by $G^\prime(\rho_i)$ for all
functions $G$  of a certain class $($namely, functions having sufficient smoothness and special properties of growth at zero and infinity$)$.}

Our main result is as follows.

{\bfseries Theorem 2.3.}  {\it Let $\Omega\subset{\mathbb R}^3$ be a bounded domain with $\partial\Omega\in C^2$, let the viscosity matrices satisfy the conditions \eqref{positiveviscos} and \eqref{triangleviscos}, let $\gamma>3$ be the adiabatic exponent,  let $\displaystyle m>\frac{2}{3}\cdot\frac{6\gamma^2-7\gamma+3}{2\gamma^2-5\gamma+1}$ be the exponent of growth of the heat-conducting coefficient, and let the other numerical parameters in the model $($$a$, $M_1$, $M_2$$)$ be arbitrary $($positive$)$. Then for all input data $\boldsymbol{f}^{(i)}\in C(\overline{\Omega})$, $i=1,2$, $\widehat{\theta}\in C^1(\partial\Omega)$ and $\widehat{\theta}>0$, Problem ${\mathcal H}$ has at least one solution in the sense of Definition 2.1.}

{\bfseries Remark 2.4.} {\it The condition imposed on $\gamma$ in Theorem 2.3 looks rather onerous and physically unrealizable. However, as experience of the development of the
Navier-Stokes-Fourier theory shows, this condition can be significantly weakened. In the paper, we regard our task as overcoming fundamental difficulties, and therefore do not dwell on this point. The same goes for the other constraints imposed on $m$ and on classes of input data.}

The rest of the paper is devoted to the proof of Theorem 2.3. The first stage is the construction of solutions of a regularized problem along with obtaining estimates (for
solutions) that are uniform with respect to the regularization parameter, and this construction is rather standard, although awkward (this is treated in sections 3 and 4). The second stage is the passage to the limit with respect to the regularization parameter using the bounds thus obtained. This is a fundamental point representing
the main difficulty, and therefore requires a more detailed exposition (see section 5).  At certain stages the proof is only sketched because, in the corresponding cases, it is
similar to the situations considered in \cite{smj12} and \cite{muchapok10}.

The regularized problem mentioned above, which we refer to as Problem ${\mathcal H}_\varepsilon$, is as follows. It is required to find functions
$\rho_i^\varepsilon$, $\boldsymbol{u}^{(i)}_\varepsilon$ and $\theta^\varepsilon$ (here and below, the superscript $\varepsilon$ is not a degree) satisfying the following equations and boundary conditions, as well as additional conditions containing the parameter $\varepsilon\in(0,1]$:
\begin{align}\label{continuityeps} -\varepsilon\Delta\rho_i^\varepsilon+{\rm div}(\rho_i^\varepsilon \boldsymbol{u}^{(i)}_\varepsilon)+\varepsilon\rho_i^\varepsilon=
 \varepsilon\frac{M_i}{|\Omega|}\quad\text{on}\quad\Omega,\quad i=1,2,\end{align}
\begin{align}\label{momentumeps}\begin{array}{c}\displaystyle
 \sum\limits_{j=1}^2 L_{ij}\boldsymbol{u}^{(j)}_\varepsilon+\frac{\varepsilon}{2}\rho_i^\varepsilon\boldsymbol{u}^{(i)}_\varepsilon+
 \frac{\varepsilon}{2}\frac{M_i}{|\Omega|}\boldsymbol{u}^{(i)}_\varepsilon+
 \frac{1}{2}\rho_i^\varepsilon(\boldsymbol{u}^{(i)}_\varepsilon\cdot\nabla)\boldsymbol{u}^{(i)}_\varepsilon+\\ \\
 \displaystyle +
 \frac{1}{2}{\rm div}(\rho_i^\varepsilon \boldsymbol{u}^{(i)}_\varepsilon\otimes \boldsymbol{u}^{(i)}_\varepsilon)+\nabla (\rho_i^\varepsilon)^\gamma+\nabla(\rho_i^\varepsilon \theta_\varepsilon)=\\ \\
 =(-1)^i a(\boldsymbol{u}^{(1)}_\varepsilon-\boldsymbol{u}^{(2)}_\varepsilon)+\rho_i^\varepsilon \boldsymbol{f}^{(i)}\quad\text{on}\quad\Omega,\quad i=1,2,\end{array}\end{align}
\begin{align}\label{energyeps}\begin{array}{c}\displaystyle -2{\rm div}\left(k(\theta^\varepsilon)\frac{\varepsilon+\theta^\varepsilon}{\theta^\varepsilon}\nabla\theta^\varepsilon\right)
 +\sum\limits_{i=1}^2\left[{\rm div}(\rho_i^\varepsilon \theta^\varepsilon\boldsymbol{u}^{(i)}_\varepsilon)+\right.\\ \\
 \displaystyle +\left.\rho_i^\varepsilon \theta^\varepsilon{\rm div}\boldsymbol{u}^{(i)}_\varepsilon
 -{\mathbb P}^{(i)}_\varepsilon:(\nabla\otimes\boldsymbol{u}^{(i)}_\varepsilon)\right]=\\ \\
 \displaystyle =a|\boldsymbol{u}^{(1)}_\varepsilon-\boldsymbol{u}^{(2)}_\varepsilon|^2+\varepsilon\gamma\sum\limits_{i=1}^2(\rho_i^\varepsilon)^{\gamma-2}|\nabla\rho_i^\varepsilon|^2
 \quad\text{on}\quad\Omega,\end{array}\end{align}
\begin{align}\label{boundvelocitydenseps}\boldsymbol{u}^{(i)}_\varepsilon=0,\quad\nabla\rho_i^\varepsilon\cdot\boldsymbol{n}=0\quad\text{on}\quad\partial\Omega,\quad i=1,2,\end{align}
\begin{align}\label{boundtempereps}2k(\theta^\varepsilon)\frac{\varepsilon+\theta^\varepsilon}{\theta^\varepsilon}\nabla\theta^\varepsilon\cdot\boldsymbol{n}+
 \varepsilon\ln\theta^\varepsilon+L(\theta^\varepsilon)(\theta^\varepsilon-\widehat{\theta})=0\quad\text{on}\quad\partial\Omega,\end{align}
\begin{align}\label{addconddensityeps}\int\limits_\Omega \rho_i^\varepsilon \,d\boldsymbol{x}=M_i,\quad i=1,2,\end{align}
where
\begin{align}\label{stresseps}{\mathbb P}^{(i)}_\varepsilon=\sum\limits_{j=1}^2\widehat{{\mathbb P}}^{(ij)}_\varepsilon,\quad
\widehat{{\mathbb P}}^{(ij)}_\varepsilon=\lambda_{ij}{\rm div}\boldsymbol{u}^{(j)}_\varepsilon{\mathbb I}+2\mu_{ij}{\mathbb D}(\boldsymbol{u}^{(j)}_\varepsilon),
\quad i,j=1,2,\end{align}
and $|\Omega|$ stands for the Lebesgue measure of $\Omega$. It can be seen that Problem ${\mathcal H}_\varepsilon$ is simply a uniformly elliptic regularization of Problem ${\mathcal H}$ with additional summands and boundary conditions imposed to preserve the useful properties of the original problem that play an important role in the theory of viscous gases, for example, the integral orthogonality of the convective terms to the velocities. The conditions \eqref{addconddensityeps} follow from \eqref{continuityeps} and \eqref{boundvelocitydenseps}.  However, we include these conditions in the formulation of Problem ${\mathcal H}_\varepsilon$ for uniformity with the original problem ${\mathcal H}$.

When studying this problem, it is sometimes convenient to use the function
\begin{align}\label{entropy}s^\varepsilon=\ln\theta^\varepsilon\end{align}
(instead of the temperature $\theta^\varepsilon$), in terms of which the relations \eqref{energyeps} and \eqref{boundtempereps} can be represented in the form
\begin{align}\label{energyentreps}\begin{array}{c}\displaystyle -2{\rm div}\left((1+e^{ms^\varepsilon})(\varepsilon+e^{s^\varepsilon})\nabla s^\varepsilon\right)=
 a|\boldsymbol{u}^{(1)}_\varepsilon-\boldsymbol{u}^{(2)}_\varepsilon|^2-\\ \\ \displaystyle
 -\sum\limits_{i=1}^2\left({\rm div}(\rho_i^\varepsilon e^{s^\varepsilon}\boldsymbol{u}^{(i)}_\varepsilon)+\rho_i^\varepsilon e^{s^\varepsilon}{\rm div}\boldsymbol{u}^{(i)}_\varepsilon
 -{\mathbb P}^{(i)}_\varepsilon:(\nabla\otimes\boldsymbol{u}^{(i)}_\varepsilon)-\right.\\ \\
 \displaystyle \left.-\varepsilon\gamma(\rho_i^\varepsilon)^{\gamma-2}|\nabla\rho_i^\varepsilon|^2\right)
 \quad\text{on}\quad\Omega,\end{array}\end{align}
\begin{align}\label{boundentreps}2(1+e^{ms^\varepsilon})(\varepsilon+e^{s^\varepsilon})\nabla s^\varepsilon\cdot\boldsymbol{n}+
 \varepsilon s^\varepsilon+L(e^{s^\varepsilon})(e^{s^\varepsilon}-\widehat{\theta})=0\quad\text{on}\quad\partial\Omega,\end{align}
respectively, and we refer to the modified problem \eqref{continuityeps}, \eqref{momentumeps}, \eqref{energyentreps}, \eqref{boundvelocitydenseps}, \eqref{boundentreps} and \eqref{addconddensityeps}
of finding functions $\rho_i^\varepsilon$, $\boldsymbol{u}^{(i)}_\varepsilon$ ($i=1,2$) and $s^\varepsilon$ as Problem $\widetilde{\mathcal H}_\varepsilon$.

\section{A priori estimates of solutions of the\\ regularized problem}

\noindent\indent We shall construct a strong solution of Problem ${\mathcal H}_\varepsilon$ in the following sense.

{\bfseries Definition 3.1.} {\it By a strong generalized solution of Problem ${\mathcal H}_\varepsilon$ we mean a pair of non-negative functions
$\rho_i^\varepsilon\in W^2_p(\Omega)$, $i=1,2,$ $p>3$, a positive function $\theta^\varepsilon\in W^2_p(\Omega)$, and a pair of vector fields
$\boldsymbol{u}^{(i)}_\varepsilon\in W^2_p(\Omega)$, $i=1,2$, that satisfy \eqref{addconddensityeps}, the equations \eqref{continuityeps}--\eqref{energyeps} almost everywhere on $\Omega$, and the conditions \eqref{boundvelocitydenseps} and \eqref{boundtempereps} almost everywhere on $\partial\Omega$.}

In the same spirit one can speak of the construction of a strong solution of Problem $\widetilde{\mathcal H}_\varepsilon$ in the following sense.

{\bfseries Definition 3.2.} {\it By a strong generalized solution of Problem $\widetilde{\mathcal H}_\varepsilon$ we mean a pair of non-negative functions
$\rho_i^\varepsilon\in W^2_p(\Omega)$, $i=1,2,$ $p>3$,  a function $s^\varepsilon\in W^2_p(\Omega)$, and a pair of vector fields
$\boldsymbol{u}^{(i)}_\varepsilon\in W^2_p(\Omega)$, $i=1,2$, satisfying the condition \eqref{addconddensityeps}, the equations \eqref{continuityeps}, \eqref{momentumeps} and \eqref{energyentreps} almost everywhere on $\Omega$, and the equations \eqref{boundvelocitydenseps} and \eqref{boundentreps} almost everywhere on $\partial\Omega$.}

Indeed, it is obvious that Definitions 3.1 and 3.2 are equivalent in view of \eqref{entropy}.

Let $C_k$, $k\in{\mathbb N}$, denote quantities depending on the objects
\begin{align}\label{criticaldata}\|\boldsymbol{f}^{(i)}\|_{C(\overline{\Omega})},\>\lambda_{ij},\>\mu_{ij},\> M_i\quad (i,j=1,2),\quad
\|\widehat{\theta}\|_{C(\partial\Omega)},\>\min\limits_{\partial\Omega}\widehat{\theta},
 \> m,\> \gamma,\> a,\quad \Omega,\end{align}
and only on these, and taking finite positive values for any data in Problem~${\mathcal H}$ that satisfy the conditions of Theorem 2.3 (thus ensuring that the objects in \eqref{criticaldata} are well defined). The quantities $C_k$ will serve as the backbone for bounds uniform with respect to $\varepsilon$  (mentioned above). It is especially important that the $C_k$ do not depend on $\varepsilon$. All the objects in \eqref{criticaldata} are numerical except for $\Omega$. We do not specify the dependence of the bounds on the geometry of the domain, although this is an interesting applied problem. If some quantity in the family $\{C_k\}$ has additional arguments, we explicitly write them out in brackets.

The object here and in section 4 is the proof of the following two statements.

{\bfseries Lemma 3.3.} {\it Under the assumptions of Theorem 2.3, for an arbitrary $\varepsilon\in(0,1]$ every strong solution of Problem  ${\mathcal H}_\varepsilon$ satisfies the inequality $($see \eqref{entropy}$)$
\begin{align}\label{mainestamate}\begin{array}{c}\displaystyle\sum\limits_{i=1}^2\left(\|\rho_i^\varepsilon\|_{L_{2\gamma}(\Omega)}+\|\boldsymbol{u}^{(i)}_\varepsilon\|_{W^1_2(\Omega)}+
\|\varepsilon\nabla\rho_i^\varepsilon\|_{L_{\frac{6\gamma}{\gamma+3}}(\Omega)}\right)+\|\theta^\varepsilon\|_{L_{3m}(\Omega)}+\\ \\ \displaystyle
+\|\nabla\theta^\varepsilon\|_{L_2(\Omega)}+\int\limits_{\partial\Omega}(e^{s^\varepsilon}+e^{-s^\varepsilon})d\sigma+
\|\nabla s^\varepsilon\|_{L_2(\Omega)}+\|\theta^\varepsilon\|_{L_{2m}(\partial\Omega)}\leqslant C_1.\end{array}\end{align}}

{\bfseries Theorem 3.4.} {\it Under the assumptions of Theorem 2.3, Problem ${\mathcal H}_\varepsilon$ has at least one strong solution for every  $\varepsilon\in(0,1]$.}

{\bfseries Proof of Lemma 3.3.}  Suppose that we have a family of functions as described in Definition 3.1 and the estimate \eqref{mainestamate} has been proved for this family. Up to the end of section 4, to avoid cumbersome formulae, we omit the superscript or subscript $\varepsilon$ in the quantities depending on this parameter, such as $\rho_i^\varepsilon$, $\boldsymbol{u}^{(i)}_\varepsilon$ ($i=1,2$), $\theta^\varepsilon$, $s^\varepsilon$, and so on.

{\itshape Stage 1:  derivation of two basic integral identities.} We take the inner product of the equations \eqref{momentumeps}  and $\boldsymbol{u}^{(i)}$,  integrate over $\Omega$, and sum over $i=1,2$.  This leads to the first of the two integral identities:
\begin{align}\label{firstidentity}\begin{array}{c}\displaystyle
\sum\limits_{i=1}^2 \int\limits_\Omega {\mathbb P}^{(i)}:(\nabla\otimes\boldsymbol{u}^{(i)})\,d\boldsymbol{x}+\frac{\varepsilon}{2}\sum\limits_{i=1}^2 \int\limits_\Omega
\rho_i |\boldsymbol{u}^{(i)}|^2 \,d\boldsymbol{x}+\\ \\
\displaystyle
+\frac{\varepsilon}{2|\Omega|}\sum\limits_{i=1}^2 M_i \int\limits_\Omega |\boldsymbol{u}^{(i)}|^2 \,d\boldsymbol{x}+\varepsilon\frac{\gamma}{\gamma-1}\sum\limits_{i=1}^2 \int\limits_\Omega \rho_i^\gamma \,d\boldsymbol{x}+\\ \\
\displaystyle +\varepsilon\gamma
\sum\limits_{i=1}^2 \int\limits_\Omega \rho_i^{\gamma-2}|\nabla\rho_i|^2 \,d\boldsymbol{x}+a\int\limits_\Omega |\boldsymbol{u}^{(1)}-\boldsymbol{u}^{(2)}|^2 \,d\boldsymbol{x}=\\ \\
\displaystyle =\frac{\varepsilon}{|\Omega|}\frac{\gamma}{\gamma-1}\sum\limits_{i=1}^2 M_i \int\limits_\Omega\rho_i^{\gamma-1}\,d\boldsymbol{x}+
\sum\limits_{i=1}^2 \int\limits_\Omega\rho_i \theta{\rm div}\boldsymbol{u}^{(i)}\,d\boldsymbol{x}+\\ \\
\displaystyle
+\sum\limits_{i=1}^2 \int\limits_\Omega \rho_i\boldsymbol{f}^{(i)}\cdot \boldsymbol{u}^{(i)}\,d\boldsymbol{x}.\end{array}\end{align}
If we multiply the equation \eqref{energyeps} by $(1-1/\theta)$, integrate over $\Omega$, and add to \eqref{firstidentity}, then we obtain the second integral identity:
\begin{align}\label{secondidentity}\begin{array}{c}\displaystyle
\sum\limits_{i=1}^2 \int\limits_\Omega\frac{{\mathbb P}^{(i)}:(\nabla\otimes\boldsymbol{u}^{(i)})}{\theta}\,d\boldsymbol{x}+2\int\limits_\Omega k(\theta)
\frac{\varepsilon+\theta}{\theta}|\nabla\ln\theta|^2 \,d\boldsymbol{x}+\\ \\
\displaystyle
+\int\limits_{\partial\Omega}L(\theta)\frac{\widehat{\theta}}{\theta}d\sigma+\int\limits_{\partial\Omega}L(\theta)\theta d\sigma+\\ \\
\displaystyle+a\int\limits_\Omega\frac{|\boldsymbol{u}^{(1)}-\boldsymbol{u}^{(2)}|^2}{\theta}\,d\boldsymbol{x}+
\frac{\varepsilon}{2}\sum\limits_{i=1}^2 \int\limits_\Omega\rho_i|\boldsymbol{u}^{(i)}|^2 \,d\boldsymbol{x}+\\ \\
\displaystyle +\frac{\varepsilon}{2|\Omega|}\sum\limits_{i=1}^2 M_i \int\limits_\Omega |\boldsymbol{u}^{(i)}|^2 \,d\boldsymbol{x}+
\varepsilon\frac{\gamma}{\gamma-1}\sum\limits_{i=1}^2 \int\limits_\Omega \rho_i^\gamma \,d\boldsymbol{x}+\\ \\
\displaystyle +\varepsilon\gamma\sum\limits_{i=1}^2 \int\limits_\Omega\frac{\rho_i^{\gamma-2}}{\theta}|\nabla\rho_i|^2 \,d\boldsymbol{x}+\varepsilon\int\limits_{\partial\Omega}(s^- e^{s^-}+s^+)d\sigma=\\ \\
\displaystyle =\sum\limits_{i=1}^2 \int\limits_\Omega(\rho_i \boldsymbol{u}^{(i)}\cdot\nabla s-
\boldsymbol{u}^{(i)}\cdot\nabla\rho_i)\,d\boldsymbol{x}+\int\limits_{\partial\Omega}L(\theta)(1+\widehat{\theta}) d\sigma+\\ \\
\displaystyle +\frac{\varepsilon}{|\Omega|}\frac{\gamma}{\gamma-1}\sum\limits_{i=1}^2 M_i \int\limits_\Omega\rho_i^{\gamma-1}\,d\boldsymbol{x}+
\sum\limits_{i=1}^2 \int\limits_\Omega \rho_i\boldsymbol{f}^{(i)}\cdot \boldsymbol{u}^{(i)}\,d\boldsymbol{x}+\\\\
\displaystyle+
\varepsilon\int\limits_{\partial\Omega}(s^+ e^{-s^+}+s^-)d\sigma,\end{array}\end{align}
where we have used the notation for the positive $z^+=z\chi(z)$ and negative  $z^-=-z\chi(-z)$  parts of an arbitrary quantity $z$ ($\chi$ stands for the Heaviside function).
We note the obvious properties  $(z\varphi(z))^+=z^+\varphi(z^+)$ and $(z\varphi(z))^-=z^-\varphi(z^-)$, which hold for an arbitrary positive function  $\varphi$.
 By \eqref{nonnegentr},  the left-hand sides of
the identities \eqref{firstidentity} and \eqref{secondidentity}  contain non-negative summands only.

{\itshape Stage 2: preliminarily estimation of the densities (using the Bogovskii operator).}
Consider functions $\boldsymbol{\varphi}^{(i)}\in \overset{\circ}{W^1_2}(\Omega)$ that are solutions of the problems
$${\rm div}\boldsymbol{\varphi}^{(i)}=\rho_i^\gamma-\frac{1}{|\Omega|}\int\limits_\Omega \rho_i^\gamma \,d\boldsymbol{x},\quad
  \boldsymbol{\varphi}^{(i)}|_{\partial\Omega}=0$$
for $i=1,2$. Under the conditions holding here, these functions exist and satisfy the estimates (see, for example, \cite{bogovskii}, \cite{novstrs04})
\begin{align}\label{bogovskii}\|\boldsymbol{\varphi}^{(i)}\|_{W^1_2(\Omega)}\leqslant C_2\|\rho_i\|_{L_{2\gamma}(\Omega)}^\gamma,\quad i=1,2.\end{align}
We use $\boldsymbol{\varphi}^{(i)}$ as test functions for \eqref{momentumeps} (that is, we take the inner product of the equation \eqref{momentumeps} and
$\boldsymbol{\varphi}^{(i)}$ and integrate over $\Omega$).  We obtain the equations
\begin{align}\label{identitybogovskii}\begin{array}{c}\displaystyle
\int\limits_\Omega \rho_i^{2\gamma} \,d\boldsymbol{x}=\frac{1}{|\Omega|}\left(\int\limits_\Omega \rho_i^\gamma \,d\boldsymbol{x}\right)^2+
\frac{1}{|\Omega|}\int\limits_\Omega \rho_i^\gamma \,d\boldsymbol{x} \int\limits_\Omega \rho_i \theta \,d\boldsymbol{x}-\\\\\displaystyle-\int\limits_\Omega \rho_i^{\gamma+1}
\theta \,d\boldsymbol{x}+\frac{\varepsilon}{2}\int\limits_\Omega \rho_i \boldsymbol{u}^{(i)}\cdot \boldsymbol{\varphi}^{(i)} \,d\boldsymbol{x}+\\ \\
\displaystyle +\varepsilon\frac{M_i}{2|\Omega|}\int\limits_\Omega \boldsymbol{u}^{(i)}\cdot \boldsymbol{\varphi}^{(i)} \,d\boldsymbol{x}+\frac{1}{2}\int\limits_\Omega (\rho_i(\boldsymbol{u}^{(i)}\cdot\nabla)\boldsymbol{u}^{(i)})\cdot
\boldsymbol{\varphi}^{(i)} \,d\boldsymbol{x}+\\ \\ \displaystyle +\sum\limits_{j=1}^2\left(\mu_{ij}\int\limits_\Omega (\nabla\otimes \boldsymbol{u}^{(j)}):
(\nabla\otimes \boldsymbol{\varphi}^{(i)})\,d\boldsymbol{x}+\right.\\\\\displaystyle+\left.(\lambda_{ij}+\mu_{ij})\int\limits_\Omega ({\rm div}\boldsymbol{u}^{(j)})({\rm div}\boldsymbol{\varphi}^{(i)})\,d\boldsymbol{x}
\right)-\\ \\ \displaystyle -\frac{1}{2}\int\limits_\Omega(\rho_i \boldsymbol{u}^{(i)}\otimes\boldsymbol{u}^{(i)}):(\nabla\otimes \boldsymbol{\varphi}^{(i)})\,d\boldsymbol{x}-\int\limits_\Omega \boldsymbol{J}^{(i)}\cdot \boldsymbol{\varphi}^{(i)}\,d\boldsymbol{x}-\\ \\
\displaystyle -\int\limits_\Omega\rho_i\boldsymbol{f}^{(i)}\cdot\boldsymbol{\varphi}^{(i)}\,d\boldsymbol{x},
\quad i=1,2.\end{array}\end{align}
Using  \eqref{addconddensityeps} and \eqref{bogovskii}  and elementary inequalities (Young's inequality with a small factor and Holder's inequality), we can estimate all the summands on the right-hand side of \eqref{identitybogovskii}  using the left-hand side and the norms of $\theta$ and~$\boldsymbol{u}^{(j)}$. This yields the estimates
\begin{align}\label{estimatebogovskii}\|\rho_i\|_{L_{2\gamma}(\Omega)}\leqslant C_3\left(1+\|\theta\|_{L_{\frac{2\gamma}{\gamma-1}}(\Omega)}^{\frac{1}{\gamma-1}}+
\sum\limits_{j=1}^2\|\boldsymbol{u}^{(j)}\|_{W^1_2(\Omega)}^{\frac{2}{\gamma-1}}\right),\quad i=1,2,\end{align}
where the condition $\gamma\geqslant 3$  has been used.

{\itshape Stage 3: a preliminary estimate for the temperature.}  To estimate the summands on the right-hand side of the identity \eqref{secondidentity}, we begin by representing the integrals in the first sum in the following form (this can easily be done using \eqref{continuityeps}):
\begin{align}\label{rhsfirstint}\begin{array}{c}\displaystyle \int\limits_\Omega(\rho_i \boldsymbol{u}^{(i)}\cdot\nabla s-
\boldsymbol{u}^{(i)}\cdot\nabla\rho_i)\,d\boldsymbol{x}=\int\limits_\Omega\left(\varepsilon\nabla\rho_i\cdot\nabla s-\varepsilon\frac{|\nabla\rho_i|^2}{\rho_i+\delta_1}\right)
\,d\boldsymbol{x}-\\ \\
\displaystyle -\varepsilon\frac{M_i}{|\Omega|}\int\limits_\Omega s \,d\boldsymbol{x}+\int\limits_\Omega \left(\varepsilon\rho_i s-\varepsilon\rho_i\ln(\rho_i+\delta_1)+\varepsilon\frac{M_i}{|\Omega|}\ln(\rho_i+\delta_1)\right)\,d\boldsymbol{x}+\\ \\
\displaystyle +\delta_1\int\limits_\Omega ({\rm div}\boldsymbol{u}^{(i)})\ln(\rho_i+\delta_1)\,d\boldsymbol{x},\end{array}\end{align}
where $\delta_1\in(0,1]$ is an arbitrary parameter. By adding the elementary inequalities
$$\varepsilon\int\limits_\Omega\nabla\rho_i\cdot\nabla s \,d\boldsymbol{x}\leqslant \frac{\varepsilon\gamma}{2}\int\limits_\Omega
\frac{(\rho_i+\delta_1)^{\gamma-2}}{\theta}|\nabla\rho_i|^2 \,d\boldsymbol{x}+$$$$+\varepsilon\int\limits_\Omega\frac{|\nabla\rho_i|^2}{\rho_i+\delta_1}\,d\boldsymbol{x}+
\frac{1}{2}\int\limits_\Omega\theta^{\frac{1}{\gamma-1}}|\nabla s|^2 \,d\boldsymbol{x},$$
$$-\varepsilon\frac{M_i}{|\Omega|}\int\limits_\Omega s \,d\boldsymbol{x}\leqslant\frac{\varepsilon}{4}\int\limits_{\partial\Omega}s^- e^{s^-}d\sigma+
\frac{1}{4}\|\nabla s\|_{L_2(\Omega)}^2+C_4,$$
$$\int\limits_\Omega \left(\varepsilon\rho_i \ln\theta-\varepsilon\rho_i\ln(\rho_i+\delta_1)+\varepsilon\frac{M_i}{|\Omega|}\ln(\rho_i+\delta_1)\right)\,d\boldsymbol{x}\leqslant
\|\rho_i \theta\|_{L_1(\Omega)}+C_5,$$
$$\delta_1\int\limits_\Omega ({\rm div}\boldsymbol{u}^{(i)})\ln(\rho_i+\delta_1)\,d\boldsymbol{x}\leqslant\int\limits_\Omega |{\rm div}\boldsymbol{u}^{(i)}|(\rho_i+1)\,d\boldsymbol{x}
\leqslant$$$$\leqslant \|\boldsymbol{u}^{(i)}\|_{W^1_2(\Omega)}^2+\frac{1}{2}\|\rho_i\|_{L_2(\Omega)}^2+\frac{|\Omega|}{2},$$
we derive a bound for the right-hand side of \eqref{rhsfirstint}. Passing to the limit as $\delta_1\to 0$ in this bound, we obtain
$$\int\limits_\Omega(\rho_i \boldsymbol{u}^{(i)}\cdot\nabla s-\boldsymbol{u}^{(i)}\cdot\nabla\rho_i)\,d\boldsymbol{x}\leqslant
  \frac{\varepsilon\gamma}{2}\int\limits_\Omega\frac{\rho_i^{\gamma-2}}{\theta}|\nabla\rho_i|^2 \,d\boldsymbol{x}+
  \frac{1}{2}\int\limits_\Omega\theta^{\frac{1}{\gamma-1}}|\nabla s|^2 \,d\boldsymbol{x}+$$
$$+\frac{\varepsilon}{4}\int\limits_{\partial\Omega}s^- e^{s^-}d\sigma+
  \frac{1}{4}\|\nabla s\|_{L_2(\Omega)}^2+\|\rho_i \theta\|_{L_1(\Omega)}+\|\boldsymbol{u}^{(i)}\|_{W^1_2(\Omega)}^2+\frac{1}{2}\|\rho_i\|_{L_2(\Omega)}^2+C_6.$$
The last integral on the right-hand side of \eqref{secondidentity} is estimated as follows:
$$\varepsilon\int\limits_{\partial\Omega}(s^+ e^{-s^+}+s^-)d\sigma\leqslant\frac{\varepsilon}{4}\int\limits_{\partial\Omega}s^- e^{s^-}d\sigma-
\varepsilon\int\limits_{\partial\Omega}s^- d\sigma+2|\partial\Omega|.$$
Using simple estimates for the other integrals, we derive (from \eqref{secondidentity}) the inequality
\begin{align}\label{estimatetempentr}\begin{array}{c}\displaystyle\int\limits_\Omega\frac{1+\theta^m}{\theta^2}|\nabla\theta|^2 \,d\boldsymbol{x}+
\int\limits_{\partial\Omega}\left(L(\theta)\theta+\frac{\widehat{\theta}}{\theta}+\varepsilon |s|\right)d\sigma\leqslant\\ \\
\displaystyle \leqslant C_7\left(\sum\limits_{i=1}^2 \|\boldsymbol{u}^{(i)}\|_{W^1_2(\Omega)}^2+\sum\limits_{i=1}^2 \|\rho_i\|_{L_2(\Omega)}^2+
\|\rho_i \theta\|_{L_1(\Omega)}+1\right).\end{array}\end{align}
This contains a bound for $\theta^{m/2}$ in $W^1_2(\Omega)$,  and thus also in  $L_6(\Omega)$, and, after further
elementary manipulations, we obtain the inequality (since $m>2$)
\begin{align}\label{estimatetemponly}
\|\theta\|_{L_{3m}(\Omega)}^m\leqslant C_8\left(\sum\limits_{i=1}^2 \|\boldsymbol{u}^{(i)}\|_{W^1_2(\Omega)}^2+\sum\limits_{i=1}^2 \|\rho_i\|_{L_2(\Omega)}^2+1\right).\end{align}

{\itshape Stage 4: estimation of the temperature and velocities using densities.} We estimate the integrals on the right-hand side of the relation \eqref{firstidentity}, but preserving only the first and fourth integrals on the left-hand side (we keep the fourth integral only temporarily, to estimate the first integral on the right-hand side). Using \eqref{coerciveviscos} and elementary inequalities, we readily obtain the bound
\begin{align}\label{estimvelocityzero}\sum\limits_{i=1}^2 \|\boldsymbol{u}^{(i)}\|_{W^1_2(\Omega)}^2\leqslant C_9\left(\sum\limits_{i=1}^2\|\rho_i\|_{L_2(\Omega)}^2+
  \sum\limits_{i=1}^2 \|\rho_i\theta\|_{L_2(\Omega)}^2+1\right),\end{align}
and, since \eqref{addconddensityeps}  can be used to estimate the norms
\begin{align}\label{estimrhotheta}\|\rho_i\theta\|_{L_2(\Omega)}^2\leqslant\delta_2\|\theta\|_{L_{3m}(\Omega)}^m+
C_{10}(\delta_2)\|\rho_i\|_{L_{2\gamma}(\Omega)}^{\frac{2\gamma}{2\gamma-1}\cdot\frac{3m+2}{3(m-2)}}\end{align}
(with an arbitrary $\delta_2>0$),  it follows from \eqref{estimvelocityzero}  that for all $\delta_3>0$
\begin{align}\label{estimatevelocityfirst}\begin{array}{c}\displaystyle\sum\limits_{i=1}^2 \|\boldsymbol{u}^{(i)}\|_{W^1_2(\Omega)}^2\leqslant \delta_3\|\theta\|_{L_{3m}(\Omega)}^m+
\\ \\ \displaystyle +C_{11}(\delta_3)\left(\sum\limits_{i=1}^2\|\rho_i\|_{L_2(\Omega)}^2+
\sum\limits_{i=1}^2\|\rho_i\|_{L_{2\gamma}(\Omega)}^{\frac{2\gamma}{2\gamma-1}\cdot\frac{3m+2}{3(m-2)}}+1\right),\end{array}\end{align}
where the first sum in the brackets can be deleted by \eqref{addconddensityeps} because $\frac{3m+2}{3(m-2)}>1$. It now follows from \eqref{estimatetemponly} that
\begin{align}\label{estimatetempfinal}\|\theta\|_{L_{3m}(\Omega)}^m\leqslant C_{12}\left(\sum\limits_{i=1}^2
\|\rho_i\|_{L_{2\gamma}(\Omega)}^{\frac{2\gamma}{2\gamma-1}\cdot\frac{3m+2}{3(m-2)}}+1\right),\end{align}
and then it also follows from \eqref{estimatevelocityfirst} that
\begin{align}\label{estimatevelosfinal}\sum\limits_{i=1}^2 \|\boldsymbol{u}^{(i)}\|_{W^1_2(\Omega)}^2\leqslant C_{13}\left(\sum\limits_{i=1}^2
\|\rho_i\|_{L_{2\gamma}(\Omega)}^{\frac{2\gamma}{2\gamma-1}\cdot\frac{3m+2}{3(m-2)}}+1\right).\end{align}

{\itshape Stage 5: completing the estimates.} By the inequality $3m>2\gamma/(2\gamma-1)$, the bound \eqref{estimatebogovskii}  can be completed using \eqref{estimatetempfinal} and
\eqref{estimatevelosfinal}:
$$\sum\limits_{i=1}^2 \|\rho_i\|_{L_{2\gamma}(\Omega)}\leqslant C_{14}\left(
\sum\limits_{i=1}^2\|\rho_i\|_{L_{2\gamma}(\Omega)}^{\frac{2\gamma}{(2\gamma-1)(\gamma-1)}\cdot\frac{3m+2}{3(m-2)}}+1\right),$$
and, since the exponent of the norm on the right-hand side is less than one, this enables us to conclude that $\displaystyle \sum\limits_{i=1}^2 \|\rho_i\|_{L_{2\gamma}(\Omega)}\leqslant C_{15}$. Then \eqref{estimatetempfinal} and \eqref{estimatevelosfinal} imply the bound
$\displaystyle \sum\limits_{i=1}^2 \|\boldsymbol{u}^{(i)}\|_{W^1_2(\Omega)}+\|\theta\|_{L_{3m}(\Omega)}\leqslant C_{16}$.
By \eqref{estimrhotheta},  the left-hand side of the inequality \eqref{estimatetempentr}  can now also be estimated, which gives, in particular, a bound for $\theta^{m/2}$ in $W^1_2(\Omega)$,  and thus in $L_4(\partial\Omega)$ as well. By elementary analysis, this shows that it only remains to justify the presence of gradients of the densities in the
estimate \eqref{mainestamate}. To this end, we invoke the following inequality \cite{lions98}, which comes from \eqref{continuityeps},  the second inequality in \eqref{boundvelocitydenseps} and \eqref{addconddensityeps}:
$$\sum\limits_{i=1}^2 \|\varepsilon\nabla\rho_i\|_{L_{\frac{6\gamma}{\gamma+3}}(\Omega)}\leqslant C_{17}\left(\sum\limits_{i=1}^2
\|\rho_i\boldsymbol{u}^{(i)}\|_{L_{\frac{6\gamma}{\gamma+3}}(\Omega)}+1\right),$$
after which the desired inequality follows from the bounds\linebreak obtained above. $\square$

\section{Solvability of the regularized problem (proof of Theorem 3.4)}

\noindent\indent In view of the observations made at the beginning of section 3, it suffices to prove the existence of a strong solution of Problem $\widetilde{\mathcal H}_\varepsilon$. This solution will be constructed as a fixed point of an operator $\mathbf{\Psi}$ formed below. Let us agree that the exponent $p>3$ is arbitrary in the present section and write $B_p(\Omega)=\{\>\boldsymbol{v}\in W^2_p(\Omega):\>\boldsymbol{v}|_{\partial\Omega}=0\>\}$.

{\itshape Stage 1: forming the main operator.} We first define several ``intermediate'' operators whose superposition will give the operator $\mathbf{\Psi}$.

We first define a pair of operators ${\mathcal R}_i$, $i=1,2$,  acting by the rule\linebreak ${\mathcal R}_i:\>\boldsymbol{w}\mapsto r$, where  $\boldsymbol{w}\in B_p(\Omega)$, and  $r$  stands for the solution of the problem
$$-\varepsilon\Delta r+{\rm div}(r\boldsymbol{w})+\varepsilon r=\varepsilon\frac{M_i}{|\Omega|},\quad \nabla r\cdot\boldsymbol{n}|_{\partial\Omega}=0.$$
Then $r={\mathcal R}_i(\boldsymbol{w})$ is non-negative \cite{novstrs04} and satisfies \eqref{addconddensityeps}.
By standard properties of elliptic boundary-value problems (see, for example, \cite{agmon1}, \cite{agmon2}), the operators
${\mathcal R}_i:\>B_p(\Omega)\to W^2_p(\Omega)$ are continuous because
$$\|{\mathcal R}_i(\boldsymbol{v})-{\mathcal R}_i(\boldsymbol{w})\|_{W^2_p(\Omega)}\leqslant
  A_1(p,\varepsilon,\|\boldsymbol{v}\|_{C^1(\overline{\Omega})},\|\boldsymbol{w}\|_{C^1(\overline{\Omega})},\Omega,M_1,M_2)\|\boldsymbol{v}-\boldsymbol{w}\|_{W^2_p(\Omega)}.$$

Another auxiliary operator is  ${\mathcal U}:\>\mathbf{g}\mapsto \mathbf{h}$, where
$\mathbf{g}=(\boldsymbol{g}^{(1)},\boldsymbol{g}^{(2)})$ has components $\boldsymbol{g}^{(i)}\in L_p(\Omega)$ and
$\mathbf{h}=(\boldsymbol{h}^{(1)},\boldsymbol{h}^{(2)})$ consists of the solutions of the problems
$$\sum\limits_{j=1}^2 L_{ij}\boldsymbol{h}^{(j)}=\boldsymbol{g}^{(i)},\quad \boldsymbol{h}^{(i)}|_{\partial\Omega}=0,\qquad i=1,2.$$
Obviously,  ${\mathcal U}:\>L_p(\Omega)\to B_p(\Omega)$ is continuous.

The third auxiliary operator is  ${\mathcal S}:\>(d,b,t)\mapsto z$, where\\
$(d,b,t)\in L_p(\Omega)\times C^1(\overline{\Omega})\times W^{1-\frac{1}{p}}_p(\partial\Omega)$ with $b>0$ and  $z$ which are the solution of the problem
$$-{\rm div}(b\nabla z)=d,\quad (b\nabla z\cdot\boldsymbol{n}+\varepsilon z)|_{\partial\Omega}=t.$$
Again by the general theory, we have the inequality
\begin{align}\label{boundopers}\|{\mathcal S}(d,b,t)\|_{W^2_p(\Omega)}\leqslant A_2(\|b\|_{C^1(\overline{\Omega})},\min\limits_{\overline{\Omega}}b,p,\Omega)
\left(\|d\|_{L_p(\Omega)}+\|t\|_{W^{1-\frac{1}{p}}_p(\partial\Omega)}\right),\end{align}
that is, ${\mathcal S}:\>L_p(\Omega)\times C^1(\overline{\Omega})\times W^{1-\frac{1}{p}}_p(\partial\Omega)\to W^2_p(\Omega)$. Applying the same bound to
the difference between two problems, that is, taking $z_k={\mathcal S}(d_k,b_k,t_k)$, $k=1,2$, and noticing that
$$(z_2-z_1)={\mathcal S}\left((d_2-d_1)+{\rm div}((b_2-b_1)\nabla z_1),b_2,(t_2-t_1)-(b_2-b_1)\frac{\partial z_1}{\partial \boldsymbol{n}}\right),$$
we obtain the following estimate from \eqref{boundopers}:
$$\|{\mathcal S}(d_2,b_2,t_2)-{\mathcal S}(d_1,b_1,t_1)\|_{W^2_p(\Omega)}\leqslant$$
$$\leqslant A_3\left(\|b_1\|_{C^1(\overline{\Omega})},\|b_2\|_{C^1(\overline{\Omega})},\min\limits_{\overline{\Omega}}b_1,\min\limits_{\overline{\Omega}}b_2,p,\Omega,
\|d_1\|_{L_p(\Omega)},\|t_1\|_{W^{1-\frac{1}{p}}_p(\partial\Omega)}\right)\times$$
$$\times\left(\|d_2-d_1\|_{L_p(\Omega)}+\|t_2-t_1\|_{W^{1-\frac{1}{p}}_p(\partial\Omega)}+\|b_2-b_1\|_{C^1(\overline{\Omega})}\right),$$
that is, the operator ${\mathcal S}$ is continuous in the same spaces.

Finally, we define the fourth family of operators ${\mathcal G}^{(1)}$, ${\mathcal G}^{(2)}$, ${\mathcal D}$, ${\mathcal B}$, and ${\mathcal T}$ as follows:
$${\mathcal G}^{(i)}(\mathbf{w},y)=-\frac{\varepsilon}{2}r_i \boldsymbol{w}^{(i)}-\frac{\varepsilon}{2}\frac{M_i}{|\Omega|}\boldsymbol{w}^{(i)}-
\frac{1}{2}r_i (\boldsymbol{w}^{(i)}\cdot\nabla)\boldsymbol{w}^{(i)}-\frac{1}{2}{\rm div}(r_i \boldsymbol{w}^{(i)}\otimes \boldsymbol{w}^{(i)})-$$
$$-\nabla r_i^\gamma-\nabla(r_i e^y)+(-1)^i a(\boldsymbol{w}^{(1)}-\boldsymbol{w}^{(2)})+r_i\boldsymbol{f}^{(i)},\quad i=1,2,$$
\smallskip
$${\mathcal D}(\mathbf{w},y)=-\sum\limits_{i=1}^2{\rm div}(r_i e^y\boldsymbol{w}^{(i)})-\sum\limits_{i=1}^2 r_i e^y{\rm div}\boldsymbol{w}^{(i)}+
a|\boldsymbol{w}^{(1)}-\boldsymbol{w}^{(2)}|^2+$$
\smallskip
$$+\sum\limits_{i=1}^2\left(\sum\limits_{j=1}^2
\left(\lambda_{ij}{\rm div}\boldsymbol{w}^{(j)}{\mathbb I}+2\mu_{ij}{\mathbb D}(\boldsymbol{w}^{(j)})\right)\right):(\nabla\otimes \boldsymbol{w}^{(i)})+
\varepsilon\gamma\sum\limits_{i=1}^2 r_i^{\gamma-2}|\nabla r_i|^2,$$
\smallskip
$${\mathcal B}(y)=2(1+e^{my})(\varepsilon+e^y),$$
$${\mathcal T}(y)=-(1+e^{(m-1)y})(e^y-\widehat{\theta})|_{\partial\Omega},$$
\smallskip
where $\mathbf{w}=(\boldsymbol{w}^{(1)},\boldsymbol{w}^{(2)})\in B_p(\Omega)$  and $r_i={\mathcal R}_i(\boldsymbol{w}^{(i)})\in W^2_p(\Omega)$, $i=1,2$,
$y\in W^2_p(\Omega)$. It can readily be seen that ${\mathcal G}^{(i)}:\>B_p(\Omega)\times W^2_p(\Omega)\to C(\overline{\Omega})$, $i=1,2$,
${\mathcal D}:\>B_p(\Omega)\times W^2_p(\Omega)\to C(\overline{\Omega})$,
${\mathcal B}:\>W^2_p(\Omega)\to C^1(\overline{\Omega})$,  and\linebreak ${\mathcal T}:\>W^2_p(\Omega)\to C^1(\partial\Omega)$. Moreover,
the following properties can easily be verified.
${\mathcal G}^{(i)}$, $i=1,2$, and ${\mathcal D}$ are well
defined, bounded, and continuous as operators from $C^1(\overline{\Omega})\times C^1(\overline{\Omega})$ to $C(\overline{\Omega})$, and
therefore these operators are compact (completely continuous) as operators from $B_p(\Omega)\times W^2_p(\Omega)$ to $L_p(\Omega)$. ${\mathcal B}$
is well defined, bounded and continuous as an operator from $C^1(\overline{\Omega})$ to $C^1(\overline{\Omega})$, and therefore ${\mathcal B}$  is compact (completely continuous) as
an operator from $W^2_p(\Omega)$ to $C^1(\overline{\Omega})$.
${\mathcal T}$  is well defined, bounded and continuous
as an operator from $C^1(\partial\Omega)$ to $C^1(\partial\Omega)$, and therefore ${\mathcal T}$  is compact (completely continuous) as an operator from
$W^2_p(\Omega)$ to $W^{1-\frac{1}{p}}_p(\partial\Omega)$, and the corresponding bounds
depend only on $p$ and the objects \eqref{criticaldata}.

We finally set $\mathbf{\Psi}=({\mathcal U}\circ({\mathcal G}^{(1)},{\mathcal G}^{(2)}),{\mathcal S}\circ({\mathcal D},{\mathcal B},{\mathcal T}))$, that is, for every
$(\mathbf{u},s)=((\boldsymbol{u}^{(1)},\boldsymbol{u}^{(2)}),s)\in B_p(\Omega)\times W^2_p(\Omega)$ we set
$$\mathbf{\Psi}((\mathbf{u},s))=\{{\mathcal U}({\mathcal G}^{(1)}(\mathbf{u},s),{\mathcal G}^{(2)}(\mathbf{u},s)),
{\mathcal S}({\mathcal D}(\mathbf{u},s),{\mathcal B}(s),{\mathcal T}(s))\}.$$
By construction, the operator  $\mathbf{\Psi}:\>B_p(\Omega)\times W^2_p(\Omega)\to B_p(\Omega)\times W^2_p(\Omega)$  is well defined,
continuous, and compact (that is, completely continuous), and the desired strong
solution of Problem $\widetilde{\mathcal H}_\varepsilon$ is of the form\linebreak
$({\mathcal R}_1(\boldsymbol{u}^{(1)}),{\mathcal R}_2(\boldsymbol{u}^{(2)}),s,\boldsymbol{u}^{(1)},\boldsymbol{u}^{(2)})$, where $(\mathbf{u},s)$ is a fixed point of $\mathbf{\Psi}$.

{\itshape Stage 2: bounds for the solutions of the operator equation.} To apply the Leray--Schauder principle \cite{gilbarg}, it remains to obtain an a priori estimate, uniform with
respect to the parameter $\lambda\in(0,1]$, for the solutions of the operator equation $\lambda\mathbf{\Psi}(\mathbf{u},s)=(\mathbf{u},s)$
 in the space $W^2_p(\Omega)$, that is, to estimate an assumed solution $(\rho_1,\rho_2,s,\boldsymbol{u}^{(1)},\boldsymbol{u}^{(2)})$ of the boundary value problem
$\widetilde{\mathcal H}_\varepsilon^{(\lambda)}$ in this space uniformly
with respect to $\lambda\in(0,1]$, where Problem $\widetilde{\mathcal H}_\varepsilon^{(\lambda)}$ consists of the relations (we omit
the subscript $\lambda$  from quantities depending  on $\lambda$)
\begin{align}\label{momentumepslam}\begin{array}{c}\displaystyle
 \sum\limits_{j=1}^2 L_{ij}\boldsymbol{u}^{(j)}+\frac{\lambda\varepsilon}{2}\rho_i\boldsymbol{u}^{(i)}+
 \frac{\lambda\varepsilon}{2}\frac{M_i}{|\Omega|}\boldsymbol{u}^{(i)}+\frac{\lambda}{2}\rho_i(\boldsymbol{u}^{(i)}\cdot\nabla)\boldsymbol{u}^{(i)}+\\\\
\displaystyle
+ \frac{\lambda}{2}{\rm div}(\rho_i \boldsymbol{u}^{(i)}\otimes \boldsymbol{u}^{(i)})+\lambda\nabla (\rho_i)^\gamma+\lambda\nabla(\rho_i e^s)=\\ \\
\displaystyle
=(-1)^i \lambda a(\boldsymbol{u}^{(1)}-\boldsymbol{u}^{(2)})+\lambda\rho_i \boldsymbol{f}^{(i)}\quad\text{on}\quad\Omega,\quad i=1,2,\end{array}\end{align}
\begin{align}\label{energyentrepslam}-2{\rm div}\left((1+e^{ms})(\varepsilon+e^{s})\nabla s\right)=\Pi\quad\text{on}\quad\Omega,\end{align}
\begin{align}\label{boundentrepslam}2(1+e^{ms})(\varepsilon+e^{s})\nabla s\cdot\boldsymbol{n}=\widehat{\Pi}\quad\text{on}\quad\partial\Omega\end{align}
together with \eqref{continuityeps}, \eqref{boundvelocitydenseps}, and \eqref{addconddensityeps}. We have used the notation
\begin{align}\label{definpi}\begin{array}{c}\displaystyle \Pi=\lambda a|\boldsymbol{u}^{(1)}-\boldsymbol{u}^{(2)}|^2-\lambda\sum\limits_{i=1}^2\left({\rm div}(\rho_i e^{s}\boldsymbol{u}^{(i)})+\right.\\ \\
 \displaystyle  \left. +\rho_i e^{s}{\rm div}\boldsymbol{u}^{(i)}-
 {\mathbb P}^{(i)}:(\nabla\otimes\boldsymbol{u}^{(i)})-\varepsilon\gamma(\rho_i)^{\gamma-2}|\nabla\rho_i|^2\right),\end{array}\end{align}
where the tensors ${\mathbb P}^{(i)}$ are defined in \eqref{stress}, and
\begin{align}\label{definpihat}\widehat{\Pi}=-\varepsilon s-\lambda L(e^{s})(e^{s}-\widehat{\theta}).\end{align}

The desired bounds are somewhat similar to those obtained in section 2,  the difference being that the parameter $\lambda$ is involved, stronger norms are to be estimated, and
the parameter $\varepsilon$  can enter dominating quantities: in this respect, the scheme of obtaining the bounds is modified. Let  $B_k$, $k\in{\mathbb N}$, denote quantities similar to $\{C_k\}$ (that is, depending on \eqref{criticaldata}) which differ from $\{C_k\}$  in two ways, namely, they can depend on  $\varepsilon>0$ and $p>3$ (the critical property is now the independence of $\lambda$). When estimating, we use both of the symbols $s$ and $\theta$, assuming that these variables are related by the formula \eqref{entropy}.

{\itshape Stage 2.1: derivation of two basic integral inequalities} (an analogue of stage 1 in the proof of Lemma 3.3). As in the derivation of the identity \eqref{firstidentity}, we take the inner product of the equations \eqref{momentumepslam} and $\boldsymbol{u}^{(i)}$, integrate over~$\Omega$, and sum over  $i=1,2$, and then we divide the result by $\lambda$. This leads to an identity coinciding with \eqref{firstidentity} with the modification that the first summand is equipped with the factor $1/\lambda\geqslant 1$. Using \eqref{nonnegentr}, we can get rid of this factor and obtain a complete analogue of \eqref{firstidentity}  with the difference that the symbol $\leqslant$ replaces the equality sign. We denote this inequality by \eqref{firstidentity}$^\prime$. As in the derivation of the identity \eqref{secondidentity}, multiply \eqref{energyentrepslam} by $(1-1/\theta)$,  integrate over~$\Omega$ after dividing by $\lambda$, and add the result to \eqref{firstidentity}$^\prime$. Again using the sign-definiteness of the integrals occurring with the factor $1/\lambda\geqslant 1$ and replacing this factor by one, we obtain a complete analogue of \eqref{secondidentity} with the difference that the symbol $\leqslant$  replaces the equality sign. We denote this inequality by~\eqref{secondidentity}$^\prime$.

{\itshape Stage 2.2: a preliminary estimate for the temperature} (an analogue of stage 3 of the proof of Lemma  3.3). Proceeding as in the derivation of the bound \eqref{estimatetempentr},  that is, estimating the expressions on the right-hand side of \eqref{secondidentity}$^\prime$ in just the same way as in Lemma 3.3, we now use the existence of integrals of  $\rho_i^\gamma$ on the left-hand side to avoid the norms $\|\rho_i\|_{L_2(\Omega)}$ on the right-hand side. This leads to the estimate
\begin{align}\label{estimatetemponlylam}\|\theta\|_{L_{3m}(\Omega)}^m\leqslant B_1\left(\sum\limits_{i=1}^2 \|\boldsymbol{u}^{(i)}\|_{W^1_2(\Omega)}^2+1\right).\end{align}

{\itshape Stage 2.3: final estimates of lower norms} (an analogue of stages 4 and 5 of the proof of Lemma 3.3). As at the previous stage, and proceeding similarly to the derivation of the bound \eqref{estimvelocityzero}, we use the existence of integrals of $\rho_i^\gamma$ on the left-hand side of \eqref{firstidentity}$^\prime$ to get rid of the norms $\|\rho_i\|_{L_2(\Omega)}$ on the right-hand side. This leads to the estimate
\begin{align}\label{estimvelocityzerolam}\sum\limits_{i=1}^2 \left(\|\boldsymbol{u}^{(i)}\|_{W^1_2(\Omega)}^2+\|\rho_i\|_{L_\gamma(\Omega)}^\gamma\right)
\leqslant B_2\left(\sum\limits_{i=1}^2 \|\rho_i\theta\|_{L_2(\Omega)}^2+1\right).\end{align}
Using the obvious inequality
$$\|\rho_i\theta\|_{L_2(\Omega)}^2\leqslant \frac{1}{2B_2}\|\rho_i\|_{L_\gamma(\Omega)}^\gamma+\frac{1}{4B_1 B_2}\|\theta\|_{L_{3m}(\Omega)}^m+B_3,$$
we can obtain (from \eqref{estimatetemponlylam} and \eqref{estimvelocityzerolam}) a bound for the left-hand side of \eqref{estimvelocityzerolam} in the form
$B_4:=4B_2 B_3+2B_2+1$,  and thus a bound for the left-hand side of \eqref{estimatetemponlylam} in the form $B_5:=B_1+B_1 B_4$. Since, as already noted above, an analogue of the
bound \eqref{estimatetempentr} holds, we obtain an analogue of the bound \eqref{mainestamate}:
\begin{align}\label{mainestamatelam}\begin{array}{c}\displaystyle\sum\limits_{i=1}^2\left(\|\rho_i\|_{L_\gamma(\Omega)}+\|\boldsymbol{u}^{(i)}\|_{W^1_2(\Omega)}\right)+
\|\theta\|_{L_{3m}(\Omega)}+\|\nabla\theta\|_{L_2(\Omega)}+\\ \\ \displaystyle
+\int\limits_{\partial\Omega}(e^s+e^{-s})d\sigma+\|\nabla s\|_{L_2(\Omega)}+\|\theta\|_{L_{2m}(\partial\Omega)}\leqslant B_6.\end{array}\end{align}

{\itshape Stage 2.4: auxiliary constructions}. We introduce the notations ($i=1,2$)
$$\boldsymbol{\alpha}^{(i)}=\frac{1}{2|\Omega|}\int\limits_\Omega \rho_i(\boldsymbol{u}^{(i)}\cdot\nabla)\boldsymbol{u}^{(i)} \,d\boldsymbol{x},$$
$$\boldsymbol{H}^{(i)}=\lambda\left(-\frac{\varepsilon}{2}\rho_i\boldsymbol{u}^{(i)}-\frac{\varepsilon M_i}{2|\Omega|}\boldsymbol{u}^{(i)}+
(-1)^i a(\boldsymbol{u}^{(1)}-\boldsymbol{u}^{(2)})+\rho_i\boldsymbol{f}^{(i)}-\boldsymbol{\alpha}^{(i)}\right),$$
$$\Phi(z)=\int\limits_0^z (1+e^{my})(\varepsilon+e^y)dy$$
and note that ${\rm sgn}\Phi(z)={\rm sgn}z$ and $|\Phi(z)|\leqslant 2+\varepsilon|z|+e^{(m+1)z}\chi(z)$. We let ${\mathbb V}^{(i)}$, $i=1,2$,  denote the solutions of the boundary-value problems
\begin{align}\label{definvi}{\rm div}{\mathbb V}^{(i)}=\frac{1}{2}\rho_i(\boldsymbol{u}^{(i)}\cdot\nabla)\boldsymbol{u}^{(i)}-\boldsymbol{\alpha}^{(i)},\quad
{\mathbb V}^{(i)}|_{\partial\Omega}=0,\end{align}
and write
$${\mathbb G}^{(i)}=\lambda\left(-\rho_i^\gamma{\mathbb I}-\rho_i \theta{\mathbb I}-\frac{1}{2}\rho_i\boldsymbol{u}^{(i)}\otimes \boldsymbol{u}^{(i)}-
{\mathbb V}^{(i)}\right),\quad i=1,2.$$
In this notation, the equations \eqref{momentumepslam}  become $\displaystyle \sum\limits_{j=1}^2 L_{ij}\boldsymbol{u}^{(j)}=\boldsymbol{H}^{(i)}+{\rm div}{\mathbb G}^{(i)}$, and the
problem \eqref{energyentrepslam}, \eqref{boundentrepslam} is represented in the form (see the notation in \eqref{definpi} and \eqref{definpihat})
\begin{align}\label{problementr} -2\Delta\Phi(s)=\Pi,\quad 2\nabla\Phi(s)\cdot\boldsymbol{n}|_{\partial\Omega}=\widehat{\Pi}.\end{align}

{\itshape Stage 2.5: estimates for the higher norms}. The route from  \eqref{mainestamatelam}  to the desired bounds for the functions $(\rho_1,\rho_2,s,\boldsymbol{u}^{(1)},\boldsymbol{u}^{(2)})$ in $W^2_p(\Omega)$ is a chain of estimates (by constants $B_k$) of the norms in the following order:
\begin{align}\label{step1}\|\boldsymbol{u}^{(i)}\|_{L_6(\Omega)},\;\|\rho_i\|_{L_\infty(\Omega)},\;\|\rho_i\boldsymbol{u}^{(i)}\|_{L_6(\Omega)},\;
             \|\boldsymbol{H}^{(i)}\|_{L_6(\Omega)},\end{align}
\begin{align}\label{step2}\|\boldsymbol{u}^{(i)}\|_{W^1_3(\Omega)},\;\|\boldsymbol{u}^{(i)}\|_{L_{4p}(\Omega)},\;\|\rho_i\|_{W^1_{4p}(\Omega)},\;
             \|\rho_i\|_{W^2_3(\Omega)},\; \|\boldsymbol{H}^{(i)}\|_{L_{4p}(\Omega)},\end{align}
\begin{align}\label{step3}\|\boldsymbol{u}^{(i)}\|_{W^1_{\min\{2p,3m\}}(\Omega)},\;\|\boldsymbol{u}^{(i)}\|_{C(\overline{\Omega})},\;
\|\rho_i\|_{W^2_{\min\{2p,3m\}}(\Omega)},\;
             \|\nabla\rho_i\|_{C^{\delta_4}(\overline{\Omega})}\;(\exists\delta_4>0),\end{align}
\begin{align}\label{step31}\|\Phi(s)\|_{W^1_2(\Omega)},\;\|\theta^{m+1}\|_{L_6(\Omega)},\;\|\theta^m \nabla\theta\|_{L_2(\Omega)},\end{align}
\begin{align}\label{step4}\|\Phi(s)\|_{W^2_2(\Omega)},\;\|s\|_{C(\overline{\Omega})},\;\|\theta\|_{C(\overline{\Omega})},\;\|\nabla s\|_{L_6(\Omega)},\;
             \|\nabla\theta\|_{L_6(\Omega)},\;\|\boldsymbol{H}^{(i)}\|_{C(\overline{\Omega})},\end{align}
\begin{align}\label{step5}\|\boldsymbol{u}^{(i)}\|_{W^2_6(\Omega)},\;\|\boldsymbol{u}^{(i)}\|_{C^{3/2}(\overline{\Omega})},
             \;\|\rho_i\|_{C^2(\overline{\Omega})},\end{align}
\begin{align}\label{step6}\|\Phi(s)\|_{W^2_6(\Omega)},\;\|\theta\|_{W^2_3(\Omega)},\;\|s\|_{W^2_3(\Omega)},\;\|\nabla\theta\|_{C(\overline{\Omega})},
             \;\|\nabla s\|_{C(\overline{\Omega})},\end{align}
\begin{align}\label{step7}\|\boldsymbol{u}^{(i)}\|_{W^2_p(\Omega)},\;\|\Phi(s)\|_{W^2_p(\Omega)},\;\|\theta\|_{W^2_p(\Omega)},\;\|s\|_{W^2_p(\Omega)}\end{align}
(with $i=1,2$ throughout), where the horizontal passages (inside each of the groups presented above) are rather trivial (they follow from the bounds obtained at the
corresponding place, embedding theorems, standard properties of elliptic problems, and properties of the function $\Phi$). We explain step by step the passages between the rows.

The bounds \eqref{mainestamatelam} and \eqref{step1} and the properties of the problem \eqref{definvi} (see \cite{bogovskii} and \cite{novstrs04})  imply the bounds for
$\|{\mathbb V}^{(i)}\|_{W^1_{3/2}(\Omega)}$, $i=1,2$,  and hence for  $\|{\mathbb G}^{(i)}\|_{L_3(\Omega)}$, $i=1,2$, which takes us to the beginning of \eqref{step2}.

Again using \eqref{definvi},  we obtain from \eqref{step2}  estimates for $\|{\mathbb V}^{(i)}\|_{L_{4p}(\Omega)}$, $i=1,2$, which readily give the beginning of \eqref{step3}.

The bound \eqref{mainestamatelam} immediately gives the bounds
\begin{align}\label{estimentrthin}\|\Phi(s)\|_{L_{\frac{2m}{m+1}}(\partial\Omega)}\leqslant B_7,\quad \int\limits_{\partial\Omega}\Phi(s)\widehat{\Pi}d\sigma\leqslant B_{8}\end{align}
and, after \eqref{step3}, the norm $\|\Pi\|_{L_2(\Omega)}$ can also be estimated.  Then \eqref{problementr}  implies an
identity whose right-hand side can be estimated by the second inequality in \eqref{estimentrthin}:
$$2\int\limits_\Omega |\nabla\Phi(s)|^2 \,d\boldsymbol{x}=\int\limits_\Omega \Phi(s)\Pi \,d\boldsymbol{x}+\int\limits_{\partial\Omega}\Phi(s)\widehat{\Pi}d\sigma\leqslant
  B_9\|\Phi(s)\|_{L_2(\Omega)}+B_8,$$
and now the first bound in \eqref{estimentrthin} enables us to pass to the beginning of \eqref{step31}.

After \eqref{step31},  the quantities $s$ and $\theta^\beta$ are estimated (for all $\beta\in[1,m+1]$) in $W^1_2(\Omega)$, and therefore in $W^{1/2}_2(\partial\Omega)$ as well. Hence,
$\|\widehat{\Pi}\|_{W^{1/2}_2(\partial\Omega)}\leqslant B_{10}$, and thus \eqref{problementr}  gives the beginning of \eqref{step4}.

The passage from \eqref{step4} to \eqref{step5}  follows from the bound\linebreak $\|{\rm div}{\mathbb G}^{(i)}\|_{L_6(\Omega)}\leqslant B_{11}$.

After \eqref{step5}, we have bounds for $\|\Pi\|_{L_6(\Omega)}$ and for $s$ and $\theta^\beta$ (for all $\beta\geqslant 1$) in $W^1_6(\Omega)$, and therefore in $W^{5/6}_6(\partial\Omega)$, and hence $\|\widehat{\Pi}\|_{W^{5/6}_6(\partial\Omega)}\leqslant B_{12}$. Thus, \eqref{problementr} gives \eqref{step6}.

The passage from \eqref{step6} to \eqref{step7}  follows from the bound\linebreak $\|{\rm div}{\mathbb G}^{(i)}\|_{C(\overline{\Omega})}\leqslant B_{13}$. If $p>6$, it is necessary to study the problem \eqref{problementr}  inside \eqref{step7}  once more. $\square$

\section{Passage to the limit as $\varepsilon\to 0$ and proof of Theorem 2.3}

\noindent\indent After using Theorem 3.4 to construct solutions $(\rho_1^\varepsilon,\rho_2^\varepsilon,\theta^\varepsilon,\boldsymbol{u}^{(1)}_\varepsilon,\boldsymbol{u}^{(2)}_\varepsilon)$ of Problems
${\mathcal H}_\varepsilon$ for all $\varepsilon\in(0,1]$  in the sense of Definition 3.1, we can apply Lemma 3.3 to these solutions, and therefore, by the bound \eqref{mainestamate}, in the above family of solutions one can choose a sequence (which we denote in the same way, that is, we do not specify the values of $\varepsilon\to 0$) in such a way that for $\varepsilon\to 0$ the following convergence relations hold (we recall the notation in \eqref{entropy}):
$$\rho_i^\varepsilon\overset{w}{\to}\rho_i\quad\text{in}\quad L_{2\gamma}(\Omega),\quad i=1,2,$$
$$\boldsymbol{u}^{(i)}_\varepsilon\overset{w}{\to}\boldsymbol{u}^{(i)}\quad\text{in}\quad W^1_2(\Omega),\quad i=1,2,$$
$$\theta^\varepsilon\overset{w}{\to}\theta\quad\text{in}\quad W^1_2(\Omega),\quad L_{3m}(\Omega)\quad \text{and}\quad L_{2m}(\partial\Omega),$$
$$s^\varepsilon\overset{w}{\to}s\quad\text{in}\quad W^1_2(\Omega),$$
\begin{align}\label{convergepressure}(\rho_i^\varepsilon)^\gamma\overset{w}{\to}\overline{\rho_i^\gamma}\quad\text{in}\quad L_2(\Omega),\quad i=1,2,\end{align}
where $((\rho_1,\rho_2),(\boldsymbol{u}^{(1)},\boldsymbol{u}^{(2)}),\theta,s,(\overline{\rho_1^\gamma},\overline{\rho_2^\gamma}))$ is an element of the space
$$(L_{2\gamma}(\Omega))^2\times (\overset{\circ}{W^1_2}(\Omega))^2\times (W^1_2(\Omega)\bigcap L_{3m}(\Omega)\bigcap L_{2m}(\partial\Omega))\times W^1_2(\Omega)\times (L_2(\Omega))^2,$$
and hence
$$\boldsymbol{u}^{(i)}_\varepsilon\to\boldsymbol{u}^{(i)}\quad\text{in}\quad L_{q_1}(\Omega),\quad\forall q_1\in[1,6),\quad i=1,2,$$
$$\theta^\varepsilon\to\theta\quad\text{in}\quad L_{q_2}(\Omega),\quad\forall q_2\in[1,3m),$$
$$\theta^\varepsilon|_{\partial\Omega}\to\theta|_{\partial\Omega}\quad\text{in}\quad L_{q_3}(\partial\Omega),\quad\forall q_3\in[1,2m),$$
$$s^\varepsilon\to s\quad\text{in}\quad L_{q_4}(\Omega),\quad\forall q_4\in[1,6);$$
here the relation $\theta=e^s$ (and thus $\theta>0$), $\rho_i\geqslant 0$, $i=1,2$, and \eqref{addconddensity} are obvious. Thus, to prove Theorem 2.3, it remains to prove the validity of the integral identities appearing in Definition 2.1. In what follows, as in \eqref{convergepressure},  we use a bar (above a symbol) to denote the weak limit of the corresponding sequence (the existence of this limit is ensured by the bounds obtained, naturally, after distinguishing a subsequence, which is immediately implied).

{\itshape Stage 1: passages to the limit, namely, the full passage in the continuity equations and a partial passage in equations for momenta}.  Multiplying \eqref{continuityeps} by $\rho_i^\varepsilon$ and integrating over $\Omega$  in view of \eqref{boundvelocitydenseps}, \eqref{mainestamate} yields the bounds  $\|\sqrt{\varepsilon}\nabla\rho_i^\varepsilon\|_{L_2(\Omega)}\leqslant C_{18}$, $i=1,2$, giving (because of the bounds for $\nabla\rho_i^\varepsilon$, contained in \eqref{mainestamate})
$$\varepsilon\nabla\rho_i^\varepsilon\to 0\quad\text{in}\quad L_{q_5}(\Omega),\quad\forall q_5\in\left[1,\frac{6\gamma}{\gamma+3}\right),\quad i=1,2.$$
Passage to the limit in the equations \eqref{continuityeps}  now becomes trivial, and we arrive
at part ${\mathcal H}1$ of Definition 2.1. Taking the inner product of the equation \eqref{momentumeps} and $\boldsymbol{\varphi}^{(i)}\in C_0^\infty(\Omega)$ and integrating over $\Omega$, we obtain weak formulations of the boundary-value problems \eqref{momentumeps}, \eqref{boundvelocitydenseps} in which one can pass to the limit, paying
heed to the above bounds and convergences, and obtain identities which differ from
those given in part ${\mathcal H}2$ of Definition 2.1 only in that the expressions $\rho_i^\gamma$ replace $\overline{\rho_i^\gamma}$ in these identities. Thus, to justify ${\mathcal H}2$, it remains to prove the equations
\begin{align}\label{limitpressure}\overline{\rho_i^\gamma}=\rho_i^\gamma,\quad i=1,2\end{align}
(which are equivalent to the strong convergence of the densities).

{\itshape Stage 2: partial passage to the limit in the energy equation}. The passage to the limit in the boundary-value problem \eqref{energyeps}, \eqref{boundtempereps}  itself, that is, in the corresponding integral identity, is not successful because of the summand
$\displaystyle\sum\limits_{i=1}^2{\mathbb P}_\varepsilon^{(i)}:(\nabla\otimes\boldsymbol{u}^{(i)}_\varepsilon)$, which is bounded uniformly with respect to $\varepsilon$ only in the space $L_1(\Omega)$ which is not weakly complete, and thus we face a problem similar to that arising in the proof of equation \eqref{limitpressure} (however, in contrast to  \eqref{limitpressure}, the present problem has no solution at the moment). Therefore, we transform the summand in question by the formula
\begin{align}\label{convertstress}\sum\limits_{i=1}^2{\mathbb P}_\varepsilon^{(i)}:(\nabla\otimes\boldsymbol{u}^{(i)}_\varepsilon)=
\sum\limits_{i=1}^2\left[{\rm div}({\mathbb P}_\varepsilon^{(i)}\boldsymbol{u}^{(i)}_\varepsilon)-\boldsymbol{u}^{(i)}_\varepsilon\cdot{\rm div}{\mathbb P}_\varepsilon^{(i)}\right],
\end{align}
in which, in turn, we express the last summand in \eqref{momentumeps} using the renormalized equations \eqref{continuityeps} (namely, the equations \eqref{continuityeps} multiplied by
$\frac{\gamma(\rho_i^\varepsilon)^{\gamma-1}}{\gamma-1}$).  In terms of integral identities, this means the following. Take an arbitrary function $\eta\in C^\infty(\overline{\Omega})$ and note that \eqref{continuityeps} and \eqref{boundvelocitydenseps}  imply the relations (for $i=1,2$)
$$0\overset{\Omega}{\simeq}{\rm div}\left(\frac{\varepsilon\gamma}{\gamma-1}(\rho_i^\varepsilon)^{\gamma-1}\eta\nabla\rho_i^\varepsilon-\frac{1}{\gamma-1}\eta
  (\rho_i^\varepsilon)^\gamma \boldsymbol{u}^{(i)}_\varepsilon\right)=$$
$$=\varepsilon\gamma\eta(\rho_i^\varepsilon)^{\gamma-2}|\nabla\rho_i^\varepsilon|^2+(\rho_i^\varepsilon)^\gamma \eta{\rm div}\boldsymbol{u}^{(i)}_\varepsilon+
  \frac{\varepsilon\gamma}{\gamma-1}(\rho_i^\varepsilon)^{\gamma-1}\nabla\rho_i^\varepsilon\cdot\nabla\eta+$$
$$+\frac{\varepsilon\gamma}{\gamma-1}\eta(\rho_i^\varepsilon)^\gamma-\frac{\varepsilon\gamma}{\gamma-1}\cdot\frac{M_i}{|\Omega|}\eta(\rho_i^\varepsilon)^{\gamma-1}-
  \frac{1}{\gamma-1}(\rho_i^\varepsilon)^\gamma\boldsymbol{u}^{(i)}_\varepsilon\cdot\nabla\eta,$$
where $\overset{\Omega}{\simeq}$  stands for coincidence up to a difference that vanishes when integrating over $\Omega$ (since this difference is the divergence of a vector field vanishing on $\partial\Omega$). Thus, we obtain the representation
\begin{align}\label{expressnablarho}\begin{array}{c} \displaystyle
\varepsilon\gamma\sum\limits_{i=1}^2 \int\limits_\Omega (\rho_i^\varepsilon)^{\gamma-2}|\nabla\rho_i^\varepsilon|^2 \eta \,d\boldsymbol{x}=
\sum\limits_{i=1}^2 \int\limits_\Omega\left(-(\rho_i^\varepsilon)^\gamma \eta{\rm div}\boldsymbol{u}^{(i)}_\varepsilon-\right.\\ \\
\displaystyle -\frac{\varepsilon\gamma}{\gamma-1}(\rho_i^\varepsilon)^{\gamma-1}\nabla\rho_i^\varepsilon\cdot\nabla\eta-
\frac{\varepsilon\gamma}{\gamma-1}\eta(\rho_i^\varepsilon)^\gamma+  \\ \\ \displaystyle \left.
+\frac{\varepsilon\gamma}{\gamma-1}\cdot\frac{M_i}{|\Omega|}\eta(\rho_i^\varepsilon)^{\gamma-1}+
  \frac{1}{\gamma-1}(\rho_i^\varepsilon)^\gamma\boldsymbol{u}^{(i)}_\varepsilon\cdot\nabla\eta\right)\,d\boldsymbol{x}.\end{array}\end{align}
We now add three integral equations:
\begin{enumerate}
  \item the equation \eqref{energyeps}, after multiplication by $\eta$ and integration over $\Omega$ paying heed to the condition \eqref{boundtempereps} (the integral formulation of \eqref{energyeps}, \eqref{boundtempereps}),
  \item  the equation \eqref{momentumeps}, after multiplication by $\boldsymbol{\varphi}^{(i)}=\eta\boldsymbol{u}^{(i)}_\varepsilon$, integration over~$\Omega$ paying heed to \eqref{boundvelocitydenseps}, and summation over $i=1,2$ (the integral representation for the last summand in \eqref{convertstress}), and
  \item  the equation \eqref{expressnablarho}, which means the use of the renormalized equations~\eqref{continuityeps}.
\end{enumerate}
This procedure gives the integral identity
$$-\sum\limits_{i=1}^2 \int\limits_\Omega \rho_i^\varepsilon\left[\frac{1}{2}|\boldsymbol{u}^{(i)}_\varepsilon|^2+\frac{1}{\gamma-1}(\rho_i^\varepsilon)^{\gamma-1}+\theta^\varepsilon
  \right]\boldsymbol{u}^{(i)}_\varepsilon\cdot\nabla\eta \,d\boldsymbol{x}-$$
$$-\sum\limits_{i=1}^2 \int\limits_\Omega \left[(\rho_i^\varepsilon)^\gamma+\rho_i^\varepsilon \theta^\varepsilon\right]\boldsymbol{u}^{(i)}_\varepsilon\cdot\nabla\eta \,
d\boldsymbol{x}+$$
$$+\sum\limits_{i=1}^2 \int\limits_\Omega {\mathbb P}_\varepsilon^{(i)}:(\boldsymbol{u}^{(i)}_\varepsilon\otimes\nabla\eta)\,d\boldsymbol{x}
=\sum\limits_{i=1}^2 \int\limits_\Omega \rho_i^\varepsilon \boldsymbol{f}^{(i)}\cdot\boldsymbol{u}^{(i)}_\varepsilon \eta \,d\boldsymbol{x}-$$
$$-2\int\limits_\Omega k(\theta^\varepsilon)\frac{\varepsilon+\theta^\varepsilon}{\theta^\varepsilon} \nabla\theta^\varepsilon\cdot\nabla\eta \,d\boldsymbol{x}-\varepsilon\int\limits_{\partial\Omega} (\ln\theta^\varepsilon)\eta\,d\sigma-$$
$$-\int\limits_{\partial\Omega}L(\theta^\varepsilon)(\theta^\varepsilon-\widehat{\theta})\eta \,d\sigma
-\frac{\varepsilon\gamma}{\gamma-1}\sum\limits_{i=1}^2 \int\limits_\Omega (\rho_i^\varepsilon)^\gamma \eta\,d\boldsymbol{x}+$$
$$+\frac{\gamma}{\gamma-1}\cdot\frac{\varepsilon}{|\Omega|}\sum\limits_{i=1}^2 M_i \int\limits_\Omega (\rho_i^\varepsilon)^{\gamma-1} \eta\,d\boldsymbol{x}
-\frac{\varepsilon}{2}\sum\limits_{i=1}^2 \int\limits_\Omega \rho_i^\varepsilon |\boldsymbol{u}^{(i)}_\varepsilon|^2 \eta\,d\boldsymbol{x}-$$
$$-\frac{\varepsilon}{2|\Omega|}\sum\limits_{i=1}^2 \int\limits_\Omega |\boldsymbol{u}^{(i)}_\varepsilon|^2 \eta\,d\boldsymbol{x}-\frac{\varepsilon\gamma}{\gamma-1}\sum\limits_{i=1}^2 \int\limits_\Omega (\rho_i^\varepsilon)^{\gamma-1}\nabla\rho_i^\varepsilon\cdot\nabla\eta\,d\boldsymbol{x},$$
which corresponds to the regularized boundary-value problem for the energy equation modified using \eqref{convertstress}. Passing to the limit in this identity, we obtain
${\mathcal H}3$ up to the relations \eqref{limitpressure}, which have still not been proved and which are thus the last obstacle in the proof of Theorem 2.3.

{\itshape Stage 3: proof of communicative relations for the effective viscous fluxes}. Consider the so-called effective viscous fluxes of the components of the mixture:
$$F_i=p_i-\sum\limits_{j=1}^2 \nu_{ij}{\rm div}\boldsymbol{u}^{(j)},\quad i=1,2,$$
the corresponding quantities for the regularized problem
\begin{align}\label{defineffviscfluxeps}F_i^\varepsilon=(\rho_i^\varepsilon)^\gamma+\rho_i^\varepsilon \theta^\varepsilon-
 \sum\limits_{j=1}^2 \nu_{ij}{\rm div}\boldsymbol{u}^{(j)}_\varepsilon,\quad i=1,2,\end{align}
and their weak limits in  $L_2(\Omega)$
$$\overline{F_i}=\overline{\rho_i^\gamma}+\rho_i \theta-\sum\limits_{j=1}^2 \nu_{ij}{\rm div}\boldsymbol{u}^{(j)},\quad i=1,2.$$
The relations \eqref{limitpressure} are equivalent to the condition that  $\overline{F_i}=F_i$, $i=1,2$ (however, in contrast to the case of densities, this is not equivalent to the strong convergence of $F_i$,  which we shall prove at this stage).

{\itshape Stage 3.1: preliminary constructions}.  We shall use the operator $\Delta^{-1}$ acting by the formula
$$(\Delta^{-1}v)(\boldsymbol{x})=-\frac{1}{4\pi}\int\limits_{{\mathbb R}^3}\frac{v(\boldsymbol{y})\,d\boldsymbol{y}}{|\boldsymbol{x}-\boldsymbol{y}|},$$
applying it to the functions $v\in L_{q_6}(\Omega)$, $q_6>3/2$, extended by zero beyond the boundary of $\Omega$. Here $\Delta^{-1}:\>L_{q_6}(\Omega)\to W^2_{q_6}(\Omega)$ and
$\Delta\circ\Delta^{-1}=I$. We also need an operator ${\rm Comm}$ acting by the formula
$${\rm Comm}(\beta,\zeta)=(\nabla\otimes\nabla\Delta^{-1}\beta)\zeta-\beta(\nabla\otimes\nabla\Delta^{-1}\zeta),$$
about which the following fact is known (see \cite{coifman}, \cite{tartar}, \cite{lions98} and \cite{feir03}).\linebreak If
$\beta_k\overset{w}{\to}0$ in $L_{q_7}(\Omega)$ and $\zeta_k\overset{w}{\to}0$ in $L_{q_8}(\Omega)$,  where
$q_7^{-1}+q_8^{-1}<1$, then ${\rm Comm}(\beta_k,\zeta_k)\overset{w}{\to}0$ in $L_{q_9}(\Omega)$, where $q_9^{-1}=q_7^{-1}+q_8^{-1}$.

For every function $\alpha\in W^2_1(\Omega)$ vanishing near $\partial\Omega$ one can readily verify the relations (see the notation in \eqref{stresseps})
$$({\rm div}\widehat{{\mathbb P}}^{(ij)}_\varepsilon)\cdot\nabla\alpha+\nu_{ij}({\rm div}\boldsymbol{u}^{(j)}_\varepsilon)\Delta\alpha\overset{\Omega}{\simeq}0,\quad i,j=1,2,$$
in which, in particular, one can take $\alpha=\tau\Delta^{-1}\omega_\varepsilon$, where $\tau\in C_0^\infty(\Omega)$ and $\omega_\varepsilon\in L_{q_6}(\Omega)$,  and finally obtain ($i,j=1,2$)
\begin{align}\label{aux1}\begin{array}{c} \displaystyle
\tau({\rm div}\widehat{{\mathbb P}}^{(ij)}_\varepsilon)\cdot\nabla\Delta^{-1}\omega_\varepsilon+\tau\nu_{ij}({\rm div}\boldsymbol{u}^{(j)}_\varepsilon)\omega_\varepsilon
\overset{\Omega}{\simeq}\\ \\
\displaystyle \overset{\Omega}{\simeq}-\nu_{ij}({\rm div}\boldsymbol{u}^{(j)}_\varepsilon)[2\nabla\tau\cdot\nabla\Delta^{-1}\omega_\varepsilon+\Delta\tau \Delta^{-1}\omega_\varepsilon]+\\\\
\displaystyle+
\widehat{{\mathbb P}}^{(ij)}_\varepsilon:[\nabla\otimes(\nabla\tau \Delta^{-1}\omega_\varepsilon)].\end{array}\end{align}
Now suppose that
\begin{align}\label{conditionomega}\omega_\varepsilon\overset{w}{\to}0\quad\text{in}\quad L_{q_6}(\Omega)\quad\text{as}\quad \varepsilon\to 0.\end{align}
Summing \eqref{aux1}  over $j=1,2$, we obtain
\begin{align}\label{aux2}\tau({\rm div}{\mathbb P}^{(i)}_\varepsilon)\cdot\nabla\Delta^{-1}\omega_\varepsilon+\tau\omega_\varepsilon\sum\limits_{j=1}^2
\nu_{ij}({\rm div}\boldsymbol{u}^{(j)}_\varepsilon)\overset{\Omega,\varepsilon}{\simeq}0,\quad i=1,2,\end{align}
where $\overset{\Omega,\varepsilon}{\simeq}$  stands for coincidence up to a difference that vanishes when integrating over $\Omega$ and passing to the limit as $\varepsilon\to 0$
(since this difference is the sum of the divergence of a vector field vanishing on $\partial\Omega$ and summands containing lower derivatives of the solution and higher derivatives of
$\tau$).

After multiplying \eqref{continuityeps} by $\tau$,  elementary manipulations yield the identities
\begin{align}\label{aux3}\begin{array}{c} \displaystyle
\nabla\Delta^{-1}{\rm div}(\tau\rho_i^\varepsilon \boldsymbol{u}^{(i)}_\varepsilon)=\tau\varepsilon\nabla\rho_i^\varepsilon+\varepsilon\nabla\Delta^{-1}
\left(\frac{\tau M_i}{|\Omega|}-\tau\rho_i^\varepsilon\right)+\\ \\ \displaystyle
+\left[\varepsilon\rho_i^\varepsilon\nabla\tau+\nabla\Delta^{-1}\left(\rho_i^\varepsilon \boldsymbol{u}^{(i)}_\varepsilon\cdot\nabla\tau-
2\varepsilon\nabla\rho_i^\varepsilon\cdot\nabla\tau-\varepsilon\rho_i^\varepsilon\Delta\tau\right)\right],\quad i=1,2.\end{array}\end{align}
Multiplying \eqref{continuityeps} by $\boldsymbol{u}^{(i)}_\varepsilon/2$ and adding to \eqref{momentumeps}, we obtain the representations
\begin{align}\label{aux4}\begin{array}{c} \displaystyle
-{\rm div}{\mathbb P}^{(i)}_\varepsilon=-\nabla [(\rho_i^\varepsilon)^\gamma+\rho_i^\varepsilon \theta_\varepsilon]
-{\rm div}(\rho_i^\varepsilon \boldsymbol{u}^{(i)}_\varepsilon\otimes \boldsymbol{u}^{(i)}_\varepsilon)+\\ \\ \displaystyle
+[(-1)^i a(\boldsymbol{u}^{(1)}_\varepsilon-\boldsymbol{u}^{(2)}_\varepsilon)+\rho_i^\varepsilon \boldsymbol{f}^{(i)}]+\frac{\varepsilon}{2}
\boldsymbol{u}^{(i)}_\varepsilon\Delta\rho_i^\varepsilon-\varepsilon\rho_i^\varepsilon\boldsymbol{u}^{(i)}_\varepsilon,\quad i=1,2.\end{array}\end{align}

{\itshape Stage 3.2: the limit of effective viscous fluxes multiplied by arbitrary functions}. We transform the expressions  $F_i^\varepsilon \omega_\varepsilon \tau$, using first the representations \eqref{defineffviscfluxeps} and then the relations \eqref{aux2} and \eqref{aux4}:
$$-F_i^\varepsilon \omega_\varepsilon \tau\overset{\Omega,\varepsilon}{\simeq}
[(-1)^i a(\boldsymbol{u}^{(1)}_\varepsilon-\boldsymbol{u}^{(2)}_\varepsilon)+\rho_i^\varepsilon \boldsymbol{f}^{(i)}]\tau\nabla\Delta^{-1}\omega_\varepsilon-
\varepsilon\tau\rho_i^\varepsilon \boldsymbol{u}^{(i)}_\varepsilon\cdot\nabla\Delta^{-1}\omega_\varepsilon+$$
$$+\boldsymbol{u}^{(i)}_\varepsilon\cdot{\rm Comm}(\omega_\varepsilon,\tau\rho_i^\varepsilon \boldsymbol{u}^{(i)}_\varepsilon)+
\boldsymbol{u}^{(i)}_\varepsilon\omega_\varepsilon\nabla\Delta^{-1}{\rm div}(\tau\rho_i^\varepsilon \boldsymbol{u}^{(i)}_\varepsilon)-$$
$$-\varepsilon(\nabla\rho_i^\varepsilon)\frac{\tau}{2}\cdot[(\nabla\otimes\boldsymbol{u}^{(i)}_\varepsilon)\nabla\Delta^{-1}\omega_\varepsilon]
-\varepsilon(\nabla\rho_i^\varepsilon)\frac{\tau}{2}\cdot[(\nabla\otimes\nabla\Delta^{-1}\omega_\varepsilon)\boldsymbol{u}^{(i)}_\varepsilon].$$
By applying  \eqref{aux3}, we finally obtain
\begin{align}\label{limiteffviscfluxabstract}\begin{array}{c} \displaystyle -F_i^\varepsilon \omega_\varepsilon \tau\overset{\Omega,\varepsilon}{\simeq}
\boldsymbol{u}^{(i)}_\varepsilon\cdot{\rm Comm}(\omega_\varepsilon,\tau\rho_i^\varepsilon \boldsymbol{u}^{(i)}_\varepsilon)+
[(-1)^i a(\boldsymbol{u}^{(1)}_\varepsilon-\boldsymbol{u}^{(2)}_\varepsilon)+\\\\\displaystyle+\rho_i^\varepsilon \boldsymbol{f}^{(i)}]\tau\nabla\Delta^{-1}\omega_\varepsilon-\left[\frac{\tau}{2}(\nabla\otimes\boldsymbol{u}^{(i)}_\varepsilon)\varepsilon\nabla\rho_i^\varepsilon+
\tau\varepsilon\rho_i^\varepsilon\boldsymbol{u}^{(i)}_\varepsilon\right]\cdot\nabla\Delta^{-1}\omega_\varepsilon-\\ \\
\displaystyle-
\frac{\tau}{2}[(\nabla\otimes\nabla\Delta^{-1}\omega_\varepsilon)\boldsymbol{u}^{(i)}_\varepsilon]\cdot\varepsilon\nabla\rho_i^\varepsilon+\\ \\
\displaystyle +\omega_\varepsilon\boldsymbol{u}^{(i)}_\varepsilon\cdot\left[\varepsilon\nabla\Delta^{-1}\left(\frac{\tau M_i}{|\Omega|}-\tau\rho_i^\varepsilon\right)+
\tau\varepsilon\nabla\rho_i^\varepsilon\right],\quad i=1,2.\end{array}\end{align}
Under the assumption that $\displaystyle q_6>\frac{6\gamma}{4\gamma-3}$,  the right-hand side of \eqref{limiteffviscfluxabstract}  converges
weakly to zero in $L_{1+\delta_5}(\Omega)$ for some $\delta_5>0$, and therefore
\begin{align}\label{limiteffviscfluxintegrabstract}\int\limits_\Omega F_i^\varepsilon \omega_\varepsilon \tau\, d\boldsymbol{x}\to 0\quad\text{as}\quad\varepsilon\to 0,\quad i=1,2.
\end{align}

{\itshape Stage 3.3: the strong convergence for effective viscous fluxes and communicative relations}. For every compact subdomain $\Omega_1\subset\subset\Omega$ we take a function~$\tau$ such that $\tau=1$ on $\Omega_1$, and $\tau\geqslant 0$ on $\Omega$. Using the relations \eqref{conditionomega} and \eqref{limiteffviscfluxintegrabstract} with $q_6=2$ and $\omega_\varepsilon=F_i^\varepsilon-\overline{F_i}$, we obtain for $i=1,2$ that
$$\int\limits_{\Omega_1}|F_i^\varepsilon-\overline{F_i}|^2 \,d\boldsymbol{x}\leqslant\int\limits_\Omega|F_i^\varepsilon-\overline{F_i}|^2 \tau\,d\boldsymbol{x}
  =$$$$=\int\limits_\Omega F_i^\varepsilon \omega_\varepsilon \tau\, d\boldsymbol{x}-\int\limits_\Omega \omega_\varepsilon \overline{F_i}\tau\, d\boldsymbol{x}\to 0\quad
  \text{as}\quad \varepsilon\to 0.$$
Due to the arbitrariness of $\Omega_1$, this means that $F_i^\varepsilon\to \overline{F_i}$ in $L_{2,{\rm loc}}(\Omega)$ and thus (since $F_i^\varepsilon$ is bounded in $L_2(\Omega)$)
also in $L_{q_{10}}(\Omega)$  for every $q_{10}<2$. This implies that if  $z_\varepsilon\overset{w}{\to}z$ in $L_{q_{11}}(\Omega)$ for some $q_{11}>2$, then
$z_\varepsilon F_i^\varepsilon\overset{w}{\to}z\overline{F_i}$ in $L_{q_{12}}(\Omega)$ for every $q_{12}<\frac{2q_{11}}{2+q_{11}}$, which yields the validity of the communicative relation $\overline{zF_i}=\overline{z}\,\overline{F_i}$. In particular, one can take $z_\varepsilon=\rho_j^\varepsilon$ with arbitrary $j=1,2$, $q_{11}=2\gamma$ and $q_{12}=1$, which implies the relations
\begin{align}\label{communicative}\begin{array}{c} \displaystyle
\int\limits_\Omega \rho_j^\varepsilon\left((\rho_i^\varepsilon)^\gamma+\rho_i^\varepsilon\theta^\varepsilon-
\sum\limits_{k=1}^2 \nu_{ik}{\rm div}\boldsymbol{u}^{(k)}_\varepsilon\right)d\boldsymbol{x}\to\\\\\displaystyle\to \int\limits_\Omega \rho_j\left(
\overline{\rho_i^\gamma}+\rho_i \theta-\sum\limits_{k=1}^2 \nu_{ik}{\rm div}\boldsymbol{u}^{(k)}\right)d\boldsymbol{x}\quad\text{as}\quad\varepsilon\to 0,\quad i,j=1,2.\end{array}\end{align}

{\itshape Stage 4: passage to the limit in the pressures} (proof of \eqref{limitpressure}, that is, the strong convergence of the densities). The relations \eqref{communicative}  are similar to a relation arising in the theory of one-component media, and the difference is not only in the number of relations (four instead of one) but also in a fundamentally new phenomenon, namely, the occurrence of mixed (dissimilar) products $\rho_j {\rm div}\boldsymbol{u}^{(k)}$, $j\neq k$,  which cannot be investigated using continuity equations (in contrast to products with $j=k$).

{\itshape Stage 4.1: renormalization and the elimination of $\rho_i {\rm div}\boldsymbol{u}^{(i)}$}. By Remark~2.2, the renormalized equations \eqref{continuity} hold. In particular, the following equations hold (in $W^{-1}_{q_{13}}(\Omega)$ for the functions $G_{\delta_6}(\rho_i)=\rho_i\ln(\rho_i+\delta_6)$ with every $\delta_6>0$ for any $q_{13}\in\left[1,\frac{6\gamma}{\gamma+3}\right)$):
$${\rm div}(G_{\delta_6}(\rho_i)\boldsymbol{u}^{(i)})+(\rho_i G^\prime_{\delta_6}(\rho_i)-G_{\delta_6}(\rho_i)){\rm div}\boldsymbol{u}^{(i)}=0,\quad i=1,2$$
(these equations are obtained by formally ``multiplying'' \eqref{continuity} by $G_{\delta_6}^\prime(\rho_i)$) and, ``integrating'' these equations over
$\Omega$ (that is, acting on the test function equal to 1) and passing to the limit as $\delta_6\to 0$, we obtain
\begin{align}\label{renormal}\int\limits_\Omega\rho_i{\rm div}\boldsymbol{u}^{(i)}\,d\boldsymbol{x}=0,\quad i=1,2.\end{align}

A similar procedure with the equations \eqref{continuityeps} (here the multiplication by $G_{\delta_6}^\prime(\rho_i^\varepsilon)$ and the integration over $\Omega$ are actually carried out and, before passage to the limit with respect to $\delta_6$,  elementary estimates are obtained) with the subsequent passage to the limit as  $\varepsilon\to 0$  leads to the inequalities
\begin{align}\label{renormallim}\int\limits_\Omega\overline{\rho_i{\rm div}\boldsymbol{u}^{(i)}}\,d\boldsymbol{x}\leqslant 0,\quad i=1,2.\end{align}
Thanks to \eqref{renormal} and \eqref{renormallim}, only dissimilar products $\rho_j {\rm div}\boldsymbol{u}^{(k)}$, $j\neq k$, actually
occur in the relations \eqref{communicative} (although the relations have the form of inequalities). In particular, for $i=j=1$  the equation \eqref{communicative} becomes
\begin{align}\label{comm11initial}\int\limits_\Omega\overline{\rho_1(\rho_1^\gamma+\rho_1\theta-\nu_{12}{\rm div}\boldsymbol{u}^{(2)})}\,d\boldsymbol{x}\leqslant
\int\limits_\Omega\rho_1(\overline{\rho_1^\gamma}+\rho_1\theta-\nu_{12}{\rm div}\boldsymbol{u}^{(2)})\,d\boldsymbol{x}.\end{align}

{\itshape Stage 4.2: proof of \eqref{limitpressure} for $i=1$}. Using \eqref{triangleviscos}, we can completely remove the velocities from \eqref{comm11initial} and hence reduce the problem to a situation similar to that in the theory of equations of a one-component medium. Thus, the subsequent actions at this stage repeat the tricks of that theory. Namely, since the function $z\mapsto z^\gamma+z\theta$ is monotone, for every $v\in L_{2\gamma}(\Omega)$, $v\geqslant 0$, we have the pointwise inequality
$(\rho_1^\varepsilon-v)((\rho_1^\varepsilon)^\gamma+\rho_1^\varepsilon \theta-v^\gamma-v\theta)\geqslant 0$, which (after integrating over $\Omega$ and passing to the limit as $\varepsilon\to 0$) takes the form
\begin{align}\label{monotone}\int\limits_\Omega\overline{\rho_1(\rho_1^\gamma+\rho_1\theta)}\,d\boldsymbol{x}\geqslant
\int\limits_\Omega v(\overline{\rho_1^\gamma}+\rho_1\theta)\,d\boldsymbol{x}+\int\limits_\Omega (\rho_1-v)(v^\gamma+v\theta)\,d\boldsymbol{x}.\end{align}
Subtracting \eqref{monotone} from \eqref{comm11initial}  and setting  $v=\rho_1+\lambda\psi$ with any $\psi\in L_{2\gamma}(\Omega)$, $\psi\geqslant 0$ and $\lambda\in{\mathbb R}^+$,
we obtain the inequality
$$\int\limits_\Omega (\overline{\rho_1^\gamma}+\rho_1\theta)\psi\,d\boldsymbol{x}\leqslant \int\limits_\Omega[(\rho_1+\lambda\psi)^\gamma+(\rho_1+\lambda\psi)\theta]
\psi\,d\boldsymbol{x}.$$
Passing to the limit as  $\lambda\to 0$  in this inequality and using the pointwise property $\overline{\rho_1^\gamma}\geqslant \rho_1^\gamma$ of weak limits \cite{ekeland}, we obtain
$(\overline{\rho_1^\gamma}-\rho_1^\gamma)\psi=0$, which gives the desired relation because $\psi$  is arbitrary. As a consequence, we see that the convergence
$\rho_1^\varepsilon\to\rho_1$ is strong in $L_\gamma(\Omega)$, and thus (by the boundedness in $L_{2\gamma}(\Omega)$) also in $L_{q_{14}}(\Omega)$ for any $q_{14}\in[1,2\gamma)$.

{\itshape Stage 4.3: the communicative relation for $\rho_2 {\rm div}\boldsymbol{u}^{(1)}$}. By repeated use of the condition \eqref{triangleviscos}, the relation \eqref{limitpressure} proved above for $i=1$, and the strong convergence of $\rho_1^\varepsilon$,  we can represent \eqref{communicative} for $i=1$ and $j=2$ in the form
\begin{align}\label{commrhodivu}\int\limits_\Omega\overline{\rho_2 {\rm div}\boldsymbol{u}^{(1)}}\,d\boldsymbol{x}=
\int\limits_\Omega\rho_2 {\rm div}\boldsymbol{u}^{(1)}\,d\boldsymbol{x}.\end{align}

{\itshape Stage 4.4: proof of \eqref{limitpressure} for $i=2$}.  Again using \eqref{renormal} and \eqref{renormallim},  we can represent \eqref{communicative} for $i=j=2$ in the form
$$\int\limits_\Omega\overline{\rho_2(\rho_2^\gamma+\rho_2\theta-\nu_{21}{\rm div}\boldsymbol{u}^{(1)})}\,d\boldsymbol{x}\leqslant
\int\limits_\Omega\rho_2(\overline{\rho_2^\gamma}+\rho_2\theta-\nu_{21}{\rm div}\boldsymbol{u}^{(1)})\,d\boldsymbol{x},$$
but this time with the velocities removed using \eqref{commrhodivu}. The rest of the argument repeats stage 4.2 verbatim but with $\rho_1$ replaced by $\rho_2$.

Thus, the relations \eqref{limitpressure} are proved, and this completes the proof of\linebreak Theorem 2.3. $\square$

\end{document}